\documentclass[11pt]{amsart}

\usepackage[T1]{fontenc}
\usepackage[utf8]{inputenc}
\usepackage{lmodern}
\usepackage{amsmath,amssymb,amsfonts,mathtools}
\usepackage[margin=1in]{geometry}
\usepackage{microtype}
\usepackage[hidelinks]{hyperref}

\allowdisplaybreaks[2]
\numberwithin{equation}{section}
\newtheorem{theorem}{Theorem}[section]
\newtheorem{proposition}[theorem]{Proposition}
\newtheorem{lemma}[theorem]{Lemma}
\newtheorem{corollary}[theorem]{Corollary}

\newcommand{\RR}{\mathbb R}
\newcommand{\OO}{\mathbb O}
\newcommand{\cH}{\mathcal H}
\newcommand{\cL}{\mathcal L}
\newcommand{\cC}{\mathcal C}
\newcommand{\F}{\mathcal F}
\newcommand{\cK}{\mathcal K}
\newcommand{\uu}{\underline u}
\newcommand{\Om}{\overline\Omega}
\DeclareMathOperator{\tr}{tr}
\DeclareMathOperator{\Sym}{Sym}
\DeclareMathOperator{\Spin}{Spin}
\DeclareMathOperator{\spin}{spin}
\DeclareMathOperator{\Real}{Re}
\newcommand{\ip}[2]{\langle #1,#2\rangle}

\title[Octonionic Dirichlet problems and Green functions]{Dirichlet problems and Green functions for the octonionic Monge--Amp\`ere equation}
\date{}

\author[Xi-Nan Ma]{Xi-Nan Ma}
\address{School of Mathematical Sciences,
University of Science and Technology of China, Hefei, China}
\email{xinan@ustc.edu.cn}

\author[Wangzhe Wu]{Wangzhe Wu}
\address{School of Mathematics and Statistics,
Ningbo University, Ningbo, Zhejiang, China}
\email{wuwz18@mail.ustc.edu.cn}

\author[Jingyong Zhu]{Jingyong Zhu}
\address{School of Mathematical Sciences, Sichuan University,
Chengdu 610065, China}
\email{jzhu@scu.edu.cn}

\begin{document}
\begin{abstract}

We study the Dirichlet problem for the octonionic Monge--Amp\`ere equation:
\begin{equation}
\begin{cases}
\det H(u)=f(q,u),\quad H(u)>0&\text{in }\Omega,\\
u=\phi&\text{on }\partial\Omega
\end{cases}
\end{equation}
where $\Omega$ is a smooth bounded domain in the octonionic plane \(\OO^2\cong\RR^{16}\).
We prove the existence of a unique smooth strictly plurisubharmonic
solution, assuming a strictly plurisubharmonic subsolution and that
\(f(q,s) >0\) is nondecreasing in \(s\). The result extends to relatively
compact domains with boundary in flat translation quotients.

On smooth bounded strictly octonionic pseudoconvex domains, we construct
the unique Green function with zero boundary values and singularity
\(-|q-a|^{-6}+O(1)\) at an interior pole \(a\). We prove local
\(C^{1,1}\) regularity up to the boundary away from the pole, together
with sharp bounds for the gradient and full real Hessian. The proof
combines boundary estimates based on translations and spin rotations
with weighted estimates obtained by differentiation along rotations
and dilations.
\end{abstract}
\maketitle

\section{Introduction}\label{sec:introduction}
Caffarelli, Nirenberg and Spruck solved the smooth Dirichlet problem for
the real Monge--Ampère equation on strictly convex domains \cite{CNS1}.
Caffarelli, Kohn, Nirenberg and Spruck proved the corresponding result
for the complex equation on strictly pseudoconvex domains \cite{CKNS}.
On a general smooth domain, a subsolution with the given boundary values
provides the barriers for the boundary estimates
\cite{GuanReal,GuanComplex,GuanCorrection}.

For the quaternionic Monge--Ampère operator introduced by Alesker
\cite{AleskerLinear,AleskerQMA}, Ma--Zhu \cite{Zhu,MaZhu}
proved a subsolution theorem on arbitrary smooth bounded domains
and on relatively compact domains in flat translation quotients,
when the right-hand side is positive and nondecreasing in the unknown function.
One application of the Ma--Zhu Theorem \cite{Zhu, MaZhu} is the periodic Dirichlet problem that
arises in OpenAI's proposed proof of Sogge's local smoothing conjecture
for the Euclidean wave equation.
In this application, a quaternionic Monge--Ampère equation in the
periodic spatial variables determines the potential used in an energy estimate
\cite[Section~1 and Theorem~A.1]{OpenAILocalSmoothing}.

This application led us to ask whether a subsolution also suffices for
the Dirichlet problem on the octonionic plane. C.~N.~Yang had suggested,
in the context of physics, that quaternions might describe symmetries of
solutions of nonlinear differential equations \cite[p.~20]{ZhangYang}.
In a later discussion of spin and internal symmetries, he also raised
the possibility that octonions might play a role \cite[p.~498]{ShiDaiYang}.

Alesker defined the octonionic Hessian and its Monge--Ampère operator
\cite{AleskerOct}. Alesker and Gordon proved an octonionic Calabi--Yau
theorem on a class of compact affine manifolds \cite{AleskerGordon},
and Wang developed the pluripotential theory on domains, including a
comparison principle \cite{Wang}. We prove that on a smooth bounded
domain in the octonionic plane, a smooth strictly plurisubharmonic
subsolution with the given boundary values suffices for smooth
solvability when the density is positive and nondecreasing in the
solution. Since octonionic multiplication is nonassociative, we use
the real action of \(\Spin(9)\) to estimate the mixed boundary derivatives.
This action
preserves the octonionic Monge--Ampère operator.
We then apply the subsolution theorem on punctured domains. Estimates
that are uniform in the density and the radius of the removed ball give
bounds for the gradient and the full real Hessian of the Green function.

\subsection{The smooth Dirichlet problem}

Let \(\OO\) denote the octonion algebra. We identify the octonionic plane
\(\OO^2:=\OO\times\OO\) with \(\RR^{16}\) as a real vector space, and
write \(q=(q_1,q_2)\).
We write 
\begin{equation}\label{eq:display-inline-011}
H(u):=(u_{\bar a b})_{a,b=1}^2
\end{equation}
for the octonionic Hessian
and use the determinant defined in Section~\ref{sec:prelim}.
A real \(C^2\) function \(u\) is called \emph{admissible} if \(H(u)>0\).

\begin{theorem}\label{thm:main}
Let \(\Omega\subset\OO^2\) be a smooth bounded domain, and let
\(\phi\in C^\infty(\partial\Omega)\).
Suppose \(f\in C^\infty(\overline\Omega\times\RR)\) satisfies
\begin{equation}\label{eq:density}
f(q,s)>0,\qquad f_s(q,s)\ge0.
\end{equation}
Assume further that there is a real function \(\underline u\in C^\infty(\overline\Omega)\)
such that
\begin{equation}\label{eq:subsolution}
H(\underline u)>0\text{ on }\overline\Omega,\qquad
\det H(\underline u)\ge f(q,\underline u)\text{ in }\Omega,\qquad
\underline u=\phi\text{ on }\partial\Omega.
\end{equation}
Then the Dirichlet problem
\begin{equation}\label{eq:dirichlet}
\begin{cases}
\det H(u)=f(q,u),\quad H(u)>0&\text{in }\Omega,\\
u=\phi&\text{on }\partial\Omega
\end{cases}
\end{equation}
has a unique admissible solution \(u\in C^\infty(\overline\Omega)\),
with \(H(u)>0\) on \(\overline\Omega\).
\end{theorem}

Equality is allowed in the determinant inequality in \eqref{eq:subsolution}.
By compactness, \(H(\underline u)\ge\kappa I_2\) on \(\overline\Omega\)
for some \(\kappa>0\). Section~\ref{sec:applications} gives examples on
annuli that admit such a subsolution but have no strictly
plurisubharmonic defining function.

The octonionic Hessian is invariant under translations, so
Theorem~\ref{thm:main} extends to flat translation quotients.

\begin{corollary}\label{cor:quotient}
Let \(X=\OO^2/\Gamma\), where \(\Gamma\subset\RR^{16}\) is a discrete
additive subgroup, and let \(M\Subset X\) have smooth boundary.
Assume that the closure of every connected component of \(M\) meets \(\partial M\).
Under \eqref{eq:density}--\eqref{eq:subsolution}, with \(M\) in place
of \(\Omega\), problem \eqref{eq:dirichlet} has a unique smooth
admissible solution on \(\overline M\).
\end{corollary}

Lemmas~\ref{prop:first-order} and \ref{lem:second-derivatives} give the
gradient bound and reduce the second derivative estimate to the boundary.
For mixed derivatives, we use affine vector fields built from translations
and \(\Spin(9)\) rotations. Their normal components vanish to second
order at the chosen boundary point, and spin invariance removes the
commutator terms from the linearized equation. The boundary barrier
then gives the mixed derivative estimate.

For the double normal derivative, we first bound the tangential
octonionic trace away from zero, using a barrier at a boundary point
where this trace is smallest. The determinant equation, together with
the tangential and mixed derivative bounds, then gives the desired
estimate. We give these boundary arguments in Section~\ref{sec:boundary}.

\subsection{The Green function}
For the homogeneous equation we prescribe zero boundary values and an
interior pole at \(a\). The theory of Harvey and Lawson gives a unique solution
continuous away from the pole on strictly pseudoconvex domains
\cite{HLdual,HLsing}. Wang identifies the radial kernel \(-|q-a|^{-6}\) and proves
local \(C^{1,1}\) regularity for the homogeneous Dirichlet problem on the
ball \cite{Wang}. We obtain bounds for the full real Hessian of the Green
function up to the boundary on strictly octonionic pseudoconvex domains.

A smooth bounded domain \(\Omega\subset\OO^2\) is called
\emph{strictly octonionic pseudoconvex} if, near each boundary point,
it admits a smooth defining function \(\rho\) with \(H(\rho)>0\).
Lemma~\ref{lem:globalrho} gives a single such defining function on
all of \(\Om\).

\begin{theorem}\label{thm:green}
Let \(\Omega\subset\OO^2\) be a smooth bounded connected strictly
octonionic pseudoconvex domain, and let \(a\in\Omega\). There is
a unique octonionic plurisubharmonic function \(G=G_\Omega(\cdot,a)\)
on \(\Omega\), continuous on \(\Om\setminus\{a\}\), such that
\begin{equation}\label{eq:greenproblem}
\begin{gathered}
G|_{\partial\Omega}=0,\qquad G(a)=-\infty,\qquad
\det H(G)=0\quad\text{in }\Omega\setminus\{a\},\\
G(q)=-|q-a|^{-6}+O(1)\quad\text{as }q\to a .
\end{gathered}
\end{equation}
Here \(\det H(G)=0\) on the closed admissible cone means that the least
eigenvalue of \(H(G)\) vanishes in the viscosity sense
(see \eqref{eq:mineigen}). Moreover, \(G\in C^{1,1}_{\mathrm{loc}}
(\Om\setminus\{a\})\), and
\begin{equation}\label{eq:mainreg}
|DG(q)|\le C|q-a|^{-7},\qquad
|D^2G(q)|\le C|q-a|^{-8}.
\end{equation}
The first estimate holds on \(\Om\setminus\{a\}\) and the second holds
almost everywhere in \(\Omega\setminus\{a\}\), with \(D^2G\) the full
real Hessian. The constant depends only on \(\Omega\), \(a\), and the
chosen strict defining function. Both powers of \(|q-a|\) are sharp.
Away from \(a\), the equation also holds almost everywhere and the
octonionic Monge--Ampère measure of \(G\) vanishes.

\end{theorem}

We normalize \(G\) by fixing the coefficient of \(-|q-a|^{-6}\).
Interpolating between the bounded remainder in \eqref{eq:greenproblem}
and the Hessian estimate in \eqref{eq:mainreg} gives the following
asymptotics at the pole.

\begin{corollary}\label{cor:profile}
Under the hypotheses of Theorem~\ref{thm:green}, let
\begin{equation}\label{eq:display-inline-059}
G_\tau(y):=\tau^6G_\Omega(a+\tau y,a).
\end{equation}

For every compact \(K\subset\RR^{16}\setminus\{0\}\),
\(0<\theta<1\), and all sufficiently small \(\tau>0\),
\begin{equation}\label{eq:profilestatement}
\bigl\|G_\tau+|\cdot|^{-6}\bigr\|_{C^{1,\theta}(K)}
\le C_{K,\theta}\tau^{3(1-\theta)} .
\end{equation}
In particular, as \(q\to a\),
\begin{equation}\label{eq:gradientstatement}
D_qG_\Omega(q,a)=6|q-a|^{-8}(q-a)+O(|q-a|^{-4}).
\end{equation}
\end{corollary}

For the complex Monge--Ampère equation, Guan
\cite{GuanComplex,GuanCorrection} and B\l ocki \cite{Blocki} studied the
regularity of the pluricomplex Green function. Our construction follows
the punctured-domain approximation of Gao, Ma and Zhang \cite{GMZ}
for homogeneous complex \(k\)-Hessian equations. For the weighted estimate
of the full real Hessian, we adapt the linear-vector-field argument of
Gao, Ma, Yu and Zhang \cite[Proposition~2.2]{GMYZ}, using the
\(\Spin(9)\) action in place of unitary rotations. To apply this argument,
we verify the invariance of the octonionic operator and choose uniformly
bounded generators for all real directions. We also obtain estimates on
both boundary components that are independent of the approximation
parameters.

To prove Theorem~\ref{thm:green}, we first solve the equation with
positive constant density on \(\Omega\setminus\overline B_r(a)\) by
Theorem~\ref{thm:main}. The lower bound for the tangential trace in
Proposition~\ref{prop:block-lower-bound} depends on the minimum of the
density. For the punctured problem we obtain a uniform lower bound on
the outer and inner boundaries. On the outer boundary this follows from
strict pseudoconvexity and the Hopf lemma. At the inner boundary we
rescale to the unit sphere. The radial subsolution bounds the normal
derivative from below and gives a fixed quadratic lower barrier near
the unit sphere. The local boundary estimates of
Section~\ref{sec:boundary} then apply on both boundaries.

For the interior estimates we differentiate the equation along
\(\Spin(9)\) rotations and dilations. The first derivative with respect to the flow parameter satisfies
the homogeneous linearized equation. Differentiating twice and using
concavity gives the differential inequality \eqref{eq:secondflowPDE}. These
vector fields span all real directions at each point away from the pole.
The maximum principle thus bounds the gradient and gives an upper bound
for every second directional derivative. Since the real Laplacian is
positive, this also bounds the full real Hessian. The bounds are uniform
as the density and the inner radius tend to zero independently.

The paper is organized as follows. Sections~\ref{sec:prelim}--\ref{sec:applications} treat the smooth
Dirichlet problem and the annular examples. We construct the Green
function in Sections~\ref{sec:construction}--\ref{sec:limit}.
Section~\ref{sec:profile} gives the approximation error
(Proposition~\ref{prop:rate}), the explicit Green function of a ball,
and the proof of Corollary~\ref{cor:profile}.

\section{Octonionic derivatives and linear algebra}\label{sec:prelim}
We follow the conventions of \cite{Baez,AleskerOct}. For real matrices,
\(|\cdot|\) denotes the Frobenius norm. On symmetric matrices we use
the associated inner product \(S:Q:=\tr(SQ)\).

\subsection{Octonions and real representations}
Write
\begin{equation}\label{eq:display-001}
\OO=\left\{\xi=\sum_{p=0}^7\xi_pe_p:\xi_p\in\RR\right\},
\qquad e_0=1,\qquad
\bar\xi=\xi_0-\sum_{p=1}^7\xi_pe_p,
\end{equation}
with \(e_1,\ldots,e_7\) an orthonormal imaginary basis.
The real part, norm, and inner product are
\begin{equation}\label{eq:display-002}
\Real\xi:=\xi_0,\qquad
|\xi|^2:=\xi\bar\xi=\sum_{p=0}^7\xi_p^2,\qquad
\ip{\xi}{\eta}:=\Real(\xi\bar\eta).
\end{equation}
We use the Cayley--Dickson representation \(\OO=\mathbb H\oplus\mathbb H\ell\), with
\begin{equation}\label{eq:display-003}
(e_0,\ldots,e_7)=(1,\mathbf i,\mathbf j,\mathbf k,\ell,
\mathbf i\ell,\mathbf j\ell,\mathbf k\ell),
\end{equation}
where \(\mathbf i,\mathbf j,\mathbf k\) are the usual quaternionic units.
For \(\xi=(\xi',\xi'')\) and \(\eta=(\eta',\eta'')\), the product is
\begin{equation}\label{eq:display-004}
(\xi',\xi'')(\eta',\eta'')
=\bigl(\xi'\eta'-\overline{\eta''}\xi'',
\eta''\xi'+\xi''\overline{\eta'}\bigr).
\end{equation}
Composition and alternativity give the following identities in \(\OO\).
\begin{equation}\label{eq:octonion-identities}
\begin{gathered}
\overline{\xi\eta}=\bar\eta\bar\xi,\qquad
|\xi\eta|=|\xi|\,|\eta|,\qquad
[\xi,\eta,\zeta]:=(\xi\eta)\zeta-\xi(\eta\zeta)
\ \text{is alternating},\\
(\xi\bar\eta)\eta=\xi|\eta|^2,\qquad
\bar\xi(\xi\eta)=|\xi|^2\eta,\qquad
\ip{\xi\eta}{\zeta}=\ip{\eta}{\bar\xi\zeta}
=\ip{\xi}{\zeta\bar\eta}.
\end{gathered}
\end{equation}
The two cancellation identities use
\(\bar\xi=2\Real\xi-\xi\) and alternativity.

Let \(R(\xi)\) denote the matrix of left multiplication by \(\xi\) in this basis, so
\begin{equation}\label{eq:display-005}
R(\xi)_{pq}:=\ip{e_p}{\xi e_q}.
\end{equation}
The inner product identities and the multiplicativity of the norm give
\begin{equation}\label{eq:left-multiplication}
R(\xi)^T=R(\bar\xi),\qquad
R(\xi)^TR(\xi)=R(\xi)R(\xi)^T=|\xi|^2I_8.
\end{equation}
For \(1\le p,q\le7\), alternativity also gives
\begin{equation}\label{eq:left-clifford}
\begin{aligned}
\bigl(R(e_p)R(e_q)+R(e_q)R(e_p)\bigr)\xi
&=(e_pe_q+e_qe_p)\xi-[e_p,e_q,\xi]-[e_q,e_p,\xi]\\
&=-2\delta_{pq}\xi.
\end{aligned}
\end{equation}
Products \(R(\xi)R(\eta)\) are ordinary real matrix products.

Let \(\operatorname{Herm}_2(\OO)\) denote the space of matrices
\begin{equation}\label{eq:determinant}
A=\begin{pmatrix}A_{11}&A_{12}\\\overline{A_{12}}&A_{22}\end{pmatrix},
\qquad
\det A:=A_{11}A_{22}-|A_{12}|^2,
\end{equation}
where \(A_{11},A_{22}\in\RR\) and \(A_{12}\in\OO\).
For \(A\in\operatorname{Herm}_2(\OO)\), define
\begin{equation}\label{eq:realization}
A_{\RR}:=
\begin{pmatrix}
A_{11}I_8&R(A_{12})\\
R(A_{12})^T&A_{22}I_8
\end{pmatrix}.
\end{equation}
We equip \(\operatorname{Herm}_2(\OO)\) with the inner product
\begin{equation}\label{eq:display-006}
\ip{A}{B}_{\OO}
:=A_{11}B_{11}+A_{22}B_{22}
+2\Real(A_{12}\overline{B_{12}}).
\end{equation}
If \(A_{11}>0\), completing the square gives
\begin{equation}\label{eq:positive-completion}
\binom{\xi}{\eta}^{\!T}A_{\RR}\binom{\xi}{\eta}
=A_{11}\left|\xi+A_{11}^{-1}R(A_{12})\eta\right|^2
+\frac{\det A}{A_{11}}|\eta|^2.
\end{equation}
It follows that
\begin{equation}\label{eq:display-007}
A>0\quad\Longleftrightarrow\quad
A_{11}>0,\ \det A>0
\quad\Longleftrightarrow\quad A_{\RR}>0.
\end{equation}
We write \(A\ge0\) when \(A\) lies in the closure of this cone.
If
\begin{equation}\label{eq:display-008}
t_A:=\frac{A_{11}+A_{22}}2,\qquad
\omega_A:=\sqrt{\left(\frac{A_{11}-A_{22}}2\right)^2+|A_{12}|^2},
\end{equation}
then 
\begin{equation}\label{eq:display-inline-089}
(A_{\RR}-t_AI_{16})^2=\omega_A^2I_{16}.
\end{equation}

The two Hermitian eigenvalues, in increasing order, satisfy
\begin{equation}\label{eq:hermitian-eigenvalues}
\begin{gathered}
\lambda_1(A):=t_A-\omega_A,\qquad \lambda_2(A):=t_A+\omega_A,\\
\tr A=\lambda_1(A)+\lambda_2(A),\qquad
\det A=\lambda_1(A)\lambda_2(A).
\end{gathered}
\end{equation}
The spectrum of \(A_{\RR}\), counted with multiplicity, consists of eight copies of each value in \eqref{eq:hermitian-eigenvalues}.

\subsection{Octonionic derivatives and the Hessian}
For \(q=(q_1,q_2)\), write
\begin{equation}\label{eq:display-009}
q_a=\sum_{p=0}^7x_{ap}e_p,\qquad a=1,2,\qquad
q_{\RR}=(x_{10},\ldots,x_{17},x_{20},\ldots,x_{27})^T.
\end{equation}
We write \(x=q_{\RR}\) when using real coordinates.
For a smooth octonion-valued function \(u\), define
\begin{equation}\label{eq:octonionic-derivatives}
u_{\bar a}:=\frac{\partial u}{\partial\bar q_a}
=\sum_{p=0}^7e_pu_{x_{ap}},\qquad
u_a:=\frac{\partial u}{\partial q_a}
=\sum_{p=0}^7u_{x_{ap}}\bar e_p.
\end{equation}
For real \(u\), \(u_a=\overline{u_{\bar a}}\), and
\begin{equation}\label{eq:hessian}
H(u)_{ab}=u_{\bar a b}
:=\frac{\partial}{\partial\bar q_a}
\left(\frac{\partial u}{\partial q_b}\right)
=\sum_{p,q=0}^7(e_p\bar e_q)u_{x_{ap}x_{bq}}
=\cH(D^2u)_{ab}.
\end{equation}
A real \(C^2\) function is plurisubharmonic if \(H(u)\ge0\) and strictly
plurisubharmonic if \(H(u)>0\).
Subscripts \(a,\bar a\) denote octonionic derivatives.
Real derivatives carry subscripts \(i,j,r=1,\ldots,16\) or \(x_{ap}\).
We sum over repeated real indices.
Since
\begin{equation}\label{eq:display-inline-102}
e_p\bar e_q+e_q\bar e_p=2\delta_{pq},
\end{equation}
symmetry of \(D^2u\) gives
\begin{equation}\label{eq:hessian-diagonal}
u_{\bar b a}=\overline{u_{\bar a b}},\qquad
u_{\bar a a}=\sum_{p=0}^7u_{x_{ap}x_{ap}},\qquad
\tr H(u)=\Delta u,\qquad H(|q|^2)=16I_2.
\end{equation}

For smooth real-valued \(g,h\) and smooth \(\chi:\RR\to\RR\), the product and chain rules give
\begin{equation}\label{eq:product-chain}
\begin{aligned}
(gh)_{\bar a b}
&=g h_{\bar a b}+h g_{\bar a b}
+g_{\bar a}h_b+h_{\bar a}g_b,\\
(\chi(g))_{\bar a b}
&=\chi'(g)g_{\bar a b}+\chi''(g)g_{\bar a}g_b.
\end{aligned}
\end{equation}
These formulas follow by differentiating
\begin{equation}\label{eq:first-product-identities}
(gh)_b=gh_b+hg_b,\qquad
(\chi(g))_b=\chi'(g)g_b,
\end{equation}
since the scalar factors are real.
For 
\begin{equation}\label{eq:display-inline-108}
\xi=(\xi_1,\xi_2)\in\OO^2\simeq\RR^{16},
\end{equation}
we have
\begin{equation}\label{eq:rank-one}
\cH(\xi\otimes\xi)=
\begin{pmatrix}
|\xi_1|^2&\xi_1\bar\xi_2\\
\xi_2\bar\xi_1&|\xi_2|^2
\end{pmatrix}\ge0,\qquad
\det\cH(\xi\otimes\xi)=0.
\end{equation}
By the spectral decomposition of a real symmetric matrix,
\(S\ge0\) implies \(\cH(S)\ge0\).

For \(S\in\Sym(16)\) and \(A,B\in\operatorname{Herm}_2(\OO)\), we have
\begin{equation}\label{eq:adjoint}
\ip{\cH(S)}{A}_{\OO}=\tr(SA_{\RR}),\qquad
\tr(A_{\RR}B_{\RR})=8\ip{A}{B}_{\OO}.
\end{equation}
The off-diagonal terms in the first identity agree because
\begin{equation}\label{eq:display-inline-113}
\ip{e_p\bar e_q}{A_{12}}=\ip{e_p}{A_{12}e_q}=R(A_{12})_{pq}.
\end{equation}

The diagonal terms agree by taking the traces of the two diagonal blocks of \(S\).
The second identity follows by polarizing \eqref{eq:left-multiplication}
and taking the trace. These identities give the normalization of the projection
\begin{equation}\label{eq:projection}
\cH(A_{\RR})=8A,\qquad
\mathsf P(S):=\frac18\cH(S)_{\RR}.
\end{equation}
The map \(\mathsf P\) is the orthogonal projection onto the space of
real representations of Hermitian octonionic matrices, since
\begin{equation}\label{eq:display-010}
\tr\bigl((S-\mathsf P(S))A_{\RR}\bigr)=0.
\end{equation}
We identify differentials with gradients using the inner product on each
space. For a differentiable function \(\Phi\) on \(\operatorname{Herm}_2(\OO)\),
the chain rule gives
\begin{equation}\label{eq:chain-rule}
D_S(\Phi\circ\cH)(S)=\bigl(D_A\Phi(\cH(S))\bigr)_{\RR}.
\end{equation}

\subsection{Mixed discriminants and the linearized operator}
Polarization of the determinant gives the mixed discriminant
\begin{equation}\label{eq:mixed-discriminant}
\begin{aligned}
\det(A,B)
&:=\frac12\bigl(A_{11}B_{22}+A_{22}B_{11}
-2\Real(A_{12}\overline{B_{12}})\bigr),\\
\det(A+tB)&=\det A+2t\det(A,B)+t^2\det B.
\end{aligned}
\end{equation}
Define the cofactor matrix by
\begin{equation}\label{eq:display-011}
A^\#:=\begin{pmatrix}A_{22}&-A_{12}\\
-\overline{A_{12}}&A_{11}\end{pmatrix}.
\end{equation}
Then
\begin{equation}\label{eq:cofactor}
2\det(A,B)=\ip{A^\#}{B}_{\OO},\qquad
A_{\RR}(A^\#)_{\RR}=(\det A)I_{16}.
\end{equation}
If \(A>0\), then \(A^\#>0\). Thus \eqref{eq:adjoint} implies
\(\det(A,B)\ge0\) for \(B\ge0\).

Set
\begin{equation}\label{eq:display-012}
\F_{\OO}(A):=\sqrt{\det A},\qquad
\F(S):=\F_{\OO}(\cH(S)),\qquad
\cC^\circ:=\{S\in\Sym(16):\cH(S)>0\}.
\end{equation}
\begin{proposition}\label{prop:structure}
The admissible cone \(\cC^\circ\) is convex.
The operator \(\F\) is smooth, elliptic, concave, and homogeneous of degree one
on \(\cC^\circ\).
For \(A=\cH(S)>0\),
\begin{equation}\label{eq:linearization}
\begin{gathered}
a=D_S\F(S)=\frac{(A^\#)_{\RR}}{2\sqrt{\det A}}>0,\qquad
D_S\F(S):Q=\frac{\det(A,\cH(Q))}{\sqrt{\det A}},\\
a:S=\F(S),\qquad
\mathcal T:=\tr a=\frac{4\tr A}{\sqrt{\det A}}\ge8.
\end{gathered}
\end{equation}
If 
\begin{equation}\label{eq:display-inline-126}
\lambda_-I_2\le A\le\lambda_+I_2
\end{equation}
for some
\(0<\lambda_-\le\lambda_+\), then
\begin{equation}\label{eq:ellipticity}
\frac12\sqrt{\frac{\lambda_-}{\lambda_+}}\,I_{16}
\le a\le
\frac12\sqrt{\frac{\lambda_+}{\lambda_-}}\,I_{16}.
\end{equation}
\end{proposition}
\begin{proof}
Write
\begin{equation}\label{eq:display-013}
A=\begin{pmatrix}t+\zeta_0&\sum_{p=0}^7\zeta_{p+1}e_p\\
\overline{\sum_{p=0}^7\zeta_{p+1}e_p}&t-\zeta_0\end{pmatrix},
\qquad \zeta=(\zeta_0,\ldots,\zeta_8).
\end{equation}
In these coordinates the positive cone is \(t>|\zeta|\), which is convex,
and 
\begin{equation}\label{eq:display-inline-129}
\F_{\OO}(A)=\sqrt{t^2-|\zeta|^2}.
\end{equation}

Let \(B\in\operatorname{Herm}_2(\OO)\) have coordinates
\((\dot t,\dot\zeta)\). Differentiating \(\F_{\OO}(A+sB)\)
twice with respect to \(s\) at \(s=0\), we obtain
\begin{equation}\label{eq:second-variation}
D^2\F_{\OO}(A)[B,B]
=\frac{\dot t^{\,2}-|\dot\zeta|^2}{\F_{\OO}(A)}
-\frac{(t\dot t-\zeta\cdot\dot\zeta)^2}{\F_{\OO}(A)^3}.
\end{equation}
For \(\zeta\ne0\), we have the following decomposition
\begin{equation}\label{eq:zeta-decomposition}
\dot\zeta=\dot\zeta_\parallel\zeta/|\zeta|+\dot\zeta_\perp,
\qquad \dot\zeta_\perp\perp\zeta,
\end{equation}
which gives us
\begin{equation}\label{eq:concavity}
\F_{\OO}(A)^3D^2\F_{\OO}(A)[B,B]
=-(|\zeta|\dot t-t\dot\zeta_\parallel)^2
-\F_{\OO}(A)^2|\dot\zeta_\perp|^2\le0.
\end{equation}
For \(\zeta=0\), the right side is \(-t^2|\dot\zeta|^2\).
Since \(\cH\) is linear, \(\F\) is concave.

From the mixed discriminant formula \eqref{eq:mixed-discriminant} and
the chain rule \eqref{eq:chain-rule}, we know 
\begin{equation}\label{eq:display-014}
D\F_{\OO}(A):B=\frac{\det(A,B)}{\sqrt{\det A}},
\qquad
a=\frac{(A^\#)_{\RR}}{2\sqrt{\det A}}.
\end{equation}
Since the eigenvalues of \(A^\#\) are \(\lambda_2(A),\lambda_1(A)\), the eigenvalues of \(a\) are
\begin{equation}\label{eq:display-015}
\frac12\sqrt{\frac{\lambda_1(A)}{\lambda_2(A)}},
\qquad
\frac12\sqrt{\frac{\lambda_2(A)}{\lambda_1(A)}},
\end{equation}
each repeated eight times. Their sum, counted with multiplicity, gives
the trace bound. The two extreme values give \eqref{eq:ellipticity}.
Homogeneity gives \(a:S=\F(S)\).
\end{proof}

Concavity and homogeneity give
\begin{equation}\label{eq:support}
\F(Z)\le D\F(S):Z,\qquad
\F_{\OO}(B)\le D\F_{\OO}(A):B.
\end{equation}
In particular,
\begin{equation}\label{eq:trace-control}
B\ge\kappa I_2
\quad\Longrightarrow\quad
D\F_{\OO}(A):B
\ge\kappa D\F_{\OO}(A):I_2=\frac{\kappa}{8}\mathcal T.
\end{equation}
At \(S=D^2u\), the linearized operator is given by
\begin{equation}\label{eq:mixed-linearization}
L\zeta=a^{ij}\zeta_{ij}
=\frac{\det(H(u),H(\zeta))}{\sqrt{\det H(u)}},
\qquad Lu=\sqrt{\det H(u)}.
\end{equation}

\subsection{Spin rotations and normal coordinates}\label{subsec:spin}
The boundary estimates use both the \(\Spin(9)\) invariance of \(\F\)
and the transitive action of \(\Spin(9)\) on \(S^{15}\). To see these
properties, consider the real symmetric matrices
\begin{equation}\label{eq:clifford-matrices}
P_0=\begin{pmatrix}I_8&0\\0&-I_8\end{pmatrix},\quad
P_1=\begin{pmatrix}0&I_8\\I_8&0\end{pmatrix},\quad
P_{p+1}=\begin{pmatrix}0&R(e_p)\\-R(e_p)&0\end{pmatrix}
\quad(1\le p\le7).
\end{equation}
By \eqref{eq:left-clifford},
\begin{equation}\label{eq:clifford-relations}
P_iP_j+P_jP_i=2\delta_{ij}I_{16},\qquad
\tr P_i=0,\qquad \tr(P_iP_j)=16\delta_{ij}.
\end{equation}
The matrices \(I_{16},P_0,\ldots,P_8\) span the range of \(\mathsf P\)
in \eqref{eq:projection}.
The skew matrices \(P_iP_j\), \(i<j\), generate the real spin
representation of \(\Spin(9)\) \cite{Baez,AleskerOct}.
The commutator identity
\begin{equation}\label{eq:display-016}
[P_iP_j,P_k]=2(\delta_{jk}P_i-\delta_{ik}P_j)
\end{equation}
shows that conjugation by \(U\in\Spin(9)\) preserves the range of
\(\mathsf P\) and its orthogonal complement. We therefore obtain
\begin{equation}\label{eq:equivariance}
\cH(USU^T)_{\RR}=U\cH(S)_{\RR}U^T,\qquad
\F(USU^T)=\F(S).
\end{equation}

To prove transitivity, we compute the action of the generators at
\(\nu_0=(0,1)\in\OO^2\):
\begin{equation}\label{eq:spin-generators}
P_0P_1\nu_0=(1,0),\qquad
P_0P_{p+1}\nu_0=(e_p,0),\qquad
P_1P_{p+1}\nu_0=(0,e_p),\quad1\le p\le7.
\end{equation}
These vectors form an orthonormal basis of
\(T_{\nu_0}S^{15}=\nu_0^\perp\).
Thus the orbit map \(U\mapsto U\nu_0\) has full rank at the identity,
and the orbit is open in \(S^{15}\).
The orbit is also closed because \(\Spin(9)\) is compact.
Since \(S^{15}\) is connected, the orbit equals \(S^{15}\).

For \(z\perp\nu_0\), we expand \(z\) in the basis
\eqref{eq:spin-generators} and use the same coefficients to form
a linear combination \(K\) of the corresponding generators.
Then \(K\nu_0=z\) and \(|K|\le C|z|\).
By transitivity, we can transfer this construction to any unit
\(\nu\) by orthogonal conjugation, which preserves the bound.
Thus
\begin{equation}\label{eq:quantitative-transitivity}
z\perp\nu\quad\Longrightarrow\quad
\exists K\in\spin(9):\quad K\nu=z,\quad |K|\le C|z|.
\end{equation}

\section{Lower-order estimates and reduction to the boundary}\label{sec:estimates}
By \eqref{eq:subsolution}, there is a constant \(\kappa>0\) such that
\begin{equation}\label{eq:kappa}
H(\underline u)\ge\kappa I_2\quad\text{on }\overline\Omega.
\end{equation}
For the continuity method we consider
\begin{equation}\label{eq:path}
\begin{gathered}
\det H(u)=F(t,q,u),\qquad u|_{\partial\Omega}=\phi,\qquad 0\le t\le1,\\
F(t,q,s):=(1-t)\det H(\underline u)(q)+tf(q,s),\qquad
g(t,q,s):=\sqrt{F(t,q,s)}.
\end{gathered}
\end{equation}
We estimate smooth admissible solutions of \(\F(D^2u)=g(t,q,u)\)
uniformly in \(t\). At \(t=0\) the function \(\underline u\) solves the
equation. It remains a subsolution for every \(t\in[0,1]\), since
\eqref{eq:density} and \eqref{eq:subsolution} give
\begin{equation}\label{eq:common-subsolution}
F_s\ge0,\qquad
\det H(\underline u)\ge F(t,q,\underline u).
\end{equation}

\subsection{Comparison and first derivatives}
\begin{lemma}[Comparison principle]\label{lem:comparison}
Let \(u_1,u_2\in C^2(\overline\Omega)\) be admissible on
\(\overline\Omega\), and let \(\widetilde g\in C^1(\overline\Omega\times\RR)\)
satisfy \(\widetilde g_s\ge0\). If
\begin{equation}\label{eq:display-017}
\F(D^2u_1)-\widetilde g(q,u_1)\ge\F(D^2u_2)-\widetilde g(q,u_2)
\quad\text{in }\Omega,\qquad
u_1\le u_2\quad\text{on }\partial\Omega,
\end{equation}
then \(u_1\le u_2\) on \(\overline\Omega\).
\end{lemma}
\begin{proof}
Set \(w:=u_1-u_2\) and 
\begin{equation}\label{eq:display-inline-182}
S_\theta:=(1-\theta)D^2u_2+\theta D^2u_1.
\end{equation}

By convexity, \(S_\theta\) lies in the admissible cone for \(0\le\theta\le1\).
Integrating the linearization along this segment gives
\begin{equation}\label{eq:averaged-linearization}
\begin{gathered}
\widehat a^{ij}w_{ij}-\widehat c\,w\ge0,\\
\widehat a^{ij}:=\int_0^1\F^{ij}(S_\theta)\,d\theta,\qquad
\widehat c:=\int_0^1\widetilde g_s(q,u_2+\theta w)\,d\theta\ge0.
\end{gathered}
\end{equation}
For fixed \(u_1,u_2\), \(S_\theta(x)\) is in a compact
subset of the admissible cone as \((x,\theta)\) ranges over
\(\overline\Omega\times[0,1]\). Thus \(\widehat a\)
is uniformly positive definite, and then the maximum principle gives \(w\le0\).
\end{proof}

Let \(h\) be the harmonic extension of \(\phi\) to \(\Omega\).
Applying Lemma~\ref{lem:comparison} to the subsolution gives the lower
bound. Since 
\begin{equation}\label{eq:display-inline-194}
\Delta(u-h)=\tr H(u)>0,
\end{equation}
 the maximum principle gives the upper bound:
\begin{equation}\label{eq:height}
\underline u\le u\le h\quad\text{on }\overline\Omega.
\end{equation}
In particular, \(u\) takes values in the fixed interval
\begin{equation}\label{eq:display-018}
I:=\left[\min_{\overline\Omega}\underline u,\max_{\partial\Omega}\phi\right].
\end{equation}
On \([0,1]\times\overline\Omega\times I\), set
\begin{equation}\label{eq:display-019}
g_-:=\min\sqrt F>0,\qquad g_+:=\max\sqrt F.
\end{equation}
Then
\begin{equation}\label{eq:path-bounds}
g_-\le g\le g_+,\qquad
0\le g_s\le C,\qquad
|D_{x,s}g|+|D_{x,s}^2g|\le C.
\end{equation}
Here \(D_{x,s}g\) and \(D_{x,s}^2g\) denote the gradient and Hessian in
the independent variables \((x,s)\). The Hessian includes \(D_x^2g\),
\(D_xg_s\), and \(g_{ss}\). Throughout the estimates we fix
\(t\in[0,1]\). Unless otherwise indicated, we evaluate \(g\) and its
derivatives at \((t,q,u(q))\), taking spatial partial derivatives with
\(s\) fixed.

For each admissible solution the linearization is elliptic on
\(\overline\Omega\). The constants \(C\) in the estimates may vary from
line to line. They depend only on \(\Omega\), the data, \(\kappa\),
\(\|\underline u\|_{C^4}\), and the bounds in \eqref{eq:path-bounds},
and are uniform in \(u\) and \(t\).

Let \(\nu\) be the inward unit normal and set \(v:=u-\underline u\).
Since \(\underline u=u=h\) on \(\partial\Omega\), \eqref{eq:height} gives
\begin{equation}\label{eq:boundary-gradient}
\underline u_\nu\le u_\nu\le h_\nu,\qquad
0\le v_\nu\le C_\nu,
\end{equation}
where \(C_\nu\ge1\) is fixed.
Since \(\phi\) prescribes the tangential derivatives of \(u\),
these bounds give
\begin{equation}\label{eq:display-020}
\|u\|_{C^0(\overline\Omega)}+\sup_{\partial\Omega}|\nabla u|\le C.
\end{equation}
Define the linearized operators at \(u\) by
\begin{equation}\label{eq:linearized-operator}
\begin{gathered}
a:=D_S\F(D^2u),\qquad \mathcal T=\tr a,\qquad c:=g_s(t,q,u)\ge0,\\
L:=a^{ij}\partial_{ij},\qquad \cL:=L-c.
\end{gathered}
\end{equation}
By \eqref{eq:mixed-linearization} and \eqref{eq:trace-control},
\begin{equation}\label{eq:basic-L}
Lu=g(t,q,u),\qquad
L\underline u\ge\frac{\kappa}{8}\mathcal T,\qquad
\mathcal T\ge8.
\end{equation}

\begin{lemma}\label{prop:first-order}
Every smooth admissible solution of \eqref{eq:path} satisfies
\begin{equation}\label{eq:first-order}
\|u\|_{C^1(\overline\Omega)}\le C.
\end{equation}
\end{lemma}
\begin{proof}
Fix a constant real unit vector \(e\) and write \(u_e:=\partial_e u\).
Differentiating the equation in the direction \(e\) gives
\begin{equation}\label{eq:first-differentiation}
L u_e=e\cdot D_xg+g_s u_e,\qquad \cL u_e=e\cdot D_xg.
\end{equation}
To apply the maximum principle, we define
\begin{equation}\label{eq:global-test}
p_0:=\underline u-1-\max_{\overline\Omega}\underline u\le-1.
\end{equation}
Since \(c\ge0\), \eqref{eq:basic-L} gives
\begin{equation}\label{eq:global-test-inequality}
\cL p_0=L\underline u-cp_0
\ge\frac{\kappa}{8}\mathcal T\ge\kappa.
\end{equation}
If 
\begin{equation}\label{eq:display-inline-226}
A_1\ge1+\sup|D_xg|/\kappa,
\end{equation}
then \eqref{eq:first-differentiation}
and \eqref{eq:global-test-inequality} give
\begin{equation}\label{eq:display-021}
\cL(\pm u_e+A_1p_0)\ge0.
\end{equation}
By the maximum principle and the boundary gradient bound, we obtain 
\begin{equation}\label{eq:display-022}
\pm u_e(x)
\le\max\{0,\sup_{\partial\Omega}(\pm u_e+A_1p_0)\}
+A_1\|p_0\|_{C^0}\le C.
\end{equation}
\end{proof}

\subsection{Reduction of second derivatives to the boundary}
\begin{lemma}\label{lem:second-derivatives}
For every smooth admissible solution of \eqref{eq:path} and every
constant real unit vector \(e\), the derivative \(u_{ee}:=\partial_e^2u\)
satisfies
\begin{equation}\label{eq:interior-C2}
\sup_\Omega u_{ee}
\le\max\{0,\sup_{\partial\Omega}u_{ee}\}+C.
\end{equation}
\end{lemma}
\begin{proof}
Differentiating \(\F(D^2u)=g(t,q,u)\) twice in the direction \(e\) gives
\begin{equation}\label{eq:display-023}
\begin{aligned}
L u_{ee}+D^2\F(D^2u)[D^2u_e,D^2u_e]
&=g_{x_ex_e}+2g_{x_es}u_e\\
&\quad+g_{ss}u_e^2+cu_{ee}.
\end{aligned}
\end{equation}
The concavity of \(\F\), the bounds \eqref{eq:path-bounds}, and the
gradient estimate \eqref{eq:first-order} therefore imply
\begin{equation}\label{eq:second-differentiation}
\begin{aligned}
\cL u_{ee}
&=g_{x_ex_e}+2g_{x_es}u_e+g_{ss}u_e^2\\
&\quad-D^2\F(D^2u)[D^2u_e,D^2u_e]\ge-C.
\end{aligned}
\end{equation}
Choose \(A_2\ge1+C/\kappa\). By \eqref{eq:global-test-inequality},
\begin{equation}\label{eq:display-024}
\cL(u_{ee}+A_2p_0)\ge-C+A_2\kappa\ge0.
\end{equation}
Since \(p_0\le0\), the maximum principle yields
\begin{equation}\label{eq:display-025}
u_{ee}(x)\le
\max\{0,\sup_{\partial\Omega}u_{ee}\}
+A_2\|p_0\|_{C^0}.
\end{equation}
\end{proof}

Since \(\Delta u>0\), an upper bound for the eigenvalues of \(D^2u\)
also gives a lower bound. By Lemma~\ref{lem:second-derivatives}, it
remains to estimate \(D^2u\) on \(\partial\Omega\).

\section{Boundary second derivative estimates}\label{sec:boundary}
In a fixed inner collar of \(\partial\Omega\), let
\begin{equation}\label{eq:display-inline-239}
d(x):=\operatorname{dist}(x,\partial\Omega).
\end{equation}
Then
\begin{equation}\label{eq:display-026}
\nu=\nabla d,\qquad |\nabla d|=1,\qquad
v=u-\underline u\ge0.
\end{equation}
In each boundary chart, a constant spin rotation sends the inward
normal at the center to \((0,1)\). The chart bounds are uniform by compactness. Since \(\mathcal T\ge8\), bounded lower-order
terms can be absorbed into \(C\mathcal T\).

\subsection{A boundary barrier}
\begin{lemma}\label{lem:collar}
There are constants \(\varepsilon_0,\eta_0,N,\delta>0\) for which
\begin{equation}\label{eq:collar-function}
V:=v+\eta_0d-\frac{N}{2}d^2
\end{equation}
satisfies
\begin{equation}\label{eq:collar-estimates}
V\ge v+\frac{\eta_0}{2}d\ge0,\qquad
\cL V\le-1-\varepsilon_0\mathcal T
\quad\text{in }\Omega_\delta:=\{0<d<\delta\}.
\end{equation}
On \(\partial\Omega\),
\begin{equation}\label{eq:collar-normal}
V=0,\qquad
\eta_0\le V_\nu=v_\nu+\eta_0\le C_\nu+\eta_0.
\end{equation}
\end{lemma}
\begin{proof}
Fix 
\begin{equation}\label{eq:display-inline-245}
C_d\ge\sup\|H(d)\|_{\rm op}
\end{equation}
in the collar and choose
\begin{equation}\label{eq:display-027}
0<\varepsilon_0\le\frac{\kappa}{32},\qquad
0<\eta_0\le\min\left\{1,\frac{\kappa}{4(1+C_d)}\right\}.
\end{equation}
Then, for \(0\le b_0\le\eta_0\),
\begin{equation}\label{eq:collar-Y}
Y:=H(\underline u)-8\varepsilon_0I_2-b_0H(d)
\ge\frac{\kappa}{2}I_2.
\end{equation}
Set 
\begin{equation}\label{eq:display-inline-247}
Z:=Y+N\cH(Dd\otimes Dd).
\end{equation}
At each point, choose a spin frame with \(Dd=(0,1)\).
Then \eqref{eq:rank-one} and \eqref{eq:equivariance} give
\begin{equation}\label{eq:collar-determinant}
Z=Y+N\begin{pmatrix}0&0\\0&1\end{pmatrix},\qquad
\det Z=\det Y+NY_{11}\ge\frac{N\kappa}{2}.
\end{equation}
Choose
\begin{equation}\label{eq:display-028}
N\ge\frac{2(g_++1)^2}{\kappa},\qquad
0<\delta\le\frac{\eta_0}{N}
\end{equation}
so that \(\Omega_\delta\) lies in the fixed collar.
With \(b_0:=\eta_0-Nd\), the product rule gives
\begin{equation}\label{eq:collar-hessian-expansion}
\begin{aligned}
H(d^2)&=2dH(d)+2\cH(Dd\otimes Dd),\\
H(V)&=H(u)-H(\underline u)+(\eta_0-Nd)H(d)
-N\cH(Dd\otimes Dd)\\
&=H(u)-Z-8\varepsilon_0I_2.
\end{aligned}
\end{equation}
Since 
\begin{equation}\label{eq:display-inline-251}
D\F_{\OO}(H(u)):I_2=\mathcal T/8,
\end{equation}
\eqref{eq:support} gives
\begin{equation}\label{eq:collar-linearization}
\begin{aligned}
LV
&=g-\frac{\det(H(u),Z)}{g}-\varepsilon_0\mathcal T\\
&\le g-\sqrt{\det Z}-\varepsilon_0\mathcal T
\le-1-\varepsilon_0\mathcal T.
\end{aligned}
\end{equation}
Since \(Nd\le\eta_0\), we have 
\begin{equation}\label{eq:display-inline-253}
V\ge v+\eta_0d/2\ge0
\end{equation}
and hence
\(\cL V\le LV\). On the boundary, \(d_\nu=1\) and
\eqref{eq:boundary-gradient} give \eqref{eq:collar-normal}.
\end{proof}

\subsection{Tangential derivatives}
For tangent vectors \(\tau,\sigma\), differentiating \(v=0\) twice
on \(\partial\Omega\) gives
\begin{equation}\label{eq:tangential-identity}
D^2v(\tau,\sigma)=v_\nu D^2d(\tau,\sigma).
\end{equation}
Indeed, along a boundary curve \(\zeta\) with \(\zeta'(0)=\tau\),
\begin{equation}\label{eq:display-029}
D^2v(\tau,\tau)+Dv\cdot\zeta''(0)=0,\qquad
D^2d(\tau,\tau)+Dd\cdot\zeta''(0)=0.
\end{equation}
Using \(Dv=v_\nu\nu\) and \(Dd=\nu\) in these identities and
polarizing, we obtain \eqref{eq:tangential-identity}. The boundary gradient estimate
\eqref{eq:boundary-gradient} now yields
\begin{equation}\label{eq:tangential-bound}
|D^2u(\tau,\sigma)|\le C|\tau|\,|\sigma|
\quad\text{on }\partial\Omega.
\end{equation}

\subsection{Mixed second derivatives}\label{sec:mixed}
Fix \(p\in\partial\Omega\) and choose \(0<r<\delta\) such that
\(B_r(p)\) lies in a boundary chart. Set
\begin{equation}\label{eq:local-comparison-domain}
D:=\Omega\cap B_r(p),\qquad d(x)\le |x-p|<\delta.
\end{equation}
Let \(b\in C^2(\overline D)\) satisfy \(b(p)=0\) and, for some
\(C>0\),
\begin{equation}\label{eq:local-comparison-assumptions}
\begin{aligned}
|b|&\le C,\quad \cL b\le C\mathcal T&&\text{in }D,\\
b(x)&\ge-C|x-p|^2&&\text{on }\partial\Omega\cap\overline D.
\end{aligned}
\end{equation}
Choose
\begin{equation}\label{eq:local-comparison-constants}
J\ge\max\{C,C/r^2\},\qquad A\varepsilon_0\ge C+2J.
\end{equation}
These choices ensure
\begin{equation}\label{eq:local-comparison-boundary}
b+J|x-p|^2+AV\ge0\quad\text{on }\partial D.
\end{equation}
Moreover, \eqref{eq:collar-estimates} and
\(\cL|x-p|^2\le2\mathcal T\) imply
\begin{equation}\label{eq:local-comparison-operator}
\cL(b+J|x-p|^2+AV)\le0\quad\text{in }D.
\end{equation}
The minimum principle gives
\begin{equation}\label{eq:local-comparison-interior}
b+J|x-p|^2+AV\ge0\quad\text{in }D.
\end{equation}
This function vanishes at \(p\), so its inward normal derivative
at \(p\) is nonnegative. By \eqref{eq:collar-normal},
\begin{equation}\label{eq:local-boundary-comparison}
b_\nu(p)\ge -A V_\nu(p)\ge -A(C_\nu+\eta_0).
\end{equation}

\begin{proposition}\label{prop:mixed}
Every smooth admissible solution of \eqref{eq:path} satisfies
\begin{equation}\label{eq:mixed-bound}
|D^2u(\tau,\nu)|\le C|\tau|
\qquad(\tau\in T\partial\Omega).
\end{equation}
\end{proposition}
\begin{proof}
Fix \(p\in\partial\Omega\) and a unit tangent vector \(\tau\) at \(p\).
The relation \(D_\tau\nu\perp\nu(p)\) and
\eqref{eq:quantitative-transitivity} allow us to choose
\(K\in\spin(9)\) with \(K\nu(p)=D_\tau\nu\) and
\(|K|\le C|D_\tau\nu|\). Use the affine field
\begin{equation}\label{eq:display-inline-291}
T(x):=\tau+K(x-p)
\end{equation}
and write \(Tv:=T\cdot Dv\).
Its normal component satisfies \(T(p)\cdot\nu(p)=0\).
For every tangent vector \(h\) at \(p\), skew symmetry of \(K\)
and symmetry of \(D\nu=D^2d\) also imply
\begin{equation}\label{eq:display-030}
D_h(T\cdot\nu)(p)
=-h\cdot K\nu(p)+\tau\cdot D_h\nu
=0.
\end{equation}
Taylor's theorem on the boundary now gives
\begin{equation}\label{eq:display-inline-298}
|T(x)\cdot\nu(x)|\le C|x-p|^2
\quad\text{on }\partial\Omega\cap B_r(p).
\end{equation}
Compactness and the bound on \(K\) make the choices of \(r\) and
\(C\) uniform in \(p\) and \(\tau\).
Using \(Dv=v_\nu\nu\) on \(\partial\Omega\), together with
\eqref{eq:boundary-gradient}, we get \(|Tv|\le C|x-p|^2\) on
\(\partial\Omega\cap B_r(p)\).

We estimate \(\cL(Tv)\) using spin invariance. Write \(S:=D^2u\). By
\eqref{eq:equivariance},
\begin{equation}\label{eq:display-inline-312}
\F(e^{sK}Se^{-sK})=\F(S).
\end{equation}
Differentiating the preceding identity at \(s=0\), we obtain
\begin{equation}\label{eq:invariance-differentiated}
0=\tr\bigl(a(KS-SK)\bigr)=-2\tr(aSK),
\end{equation}
where the last equality uses \(a^T=a\), \(S^T=S\) and \(K^T=-K\).
The affine field satisfies \(\partial_iT^j=K_{ji}\) and \(D^2T=0\).
The product rule and \eqref{eq:first-differentiation}, together with
\eqref{eq:invariance-differentiated}, therefore give
\begin{equation}\label{eq:display-031}
\cL(Tu)=T\cdot D_xg+2\sum_{i,j}K_{ji}(aS)_{ij}
=T\cdot D_xg.
\end{equation}
The bounded derivatives of \(T\) and \(\uu\) give
\(|\cL(T\uu)|\le C\mathcal T\), so \(|\cL(Tv)|\le C\mathcal T\).
We also have \(|Tv|\le C\) by \eqref{eq:first-order}.

Since \(Tv(p)=0\), applying \eqref{eq:local-boundary-comparison}
to \(b=\pm Tv\) gives \(|\partial_\nu(Tv)(p)|\le C\).
At \(p\),
\begin{equation}\label{eq:display-032}
\partial_\nu(Tv)=D^2v(\nu,\tau)+v_\nu(K\nu)\cdot\nu
=D^2v(\nu,\tau).
\end{equation}
Adding \(D^2\uu(\nu,\tau)\) proves \eqref{eq:mixed-bound}.
\end{proof}

\subsection{The double normal derivative}\label{sec:normal}
The determinant equation bounds the double normal derivative once the
tangential octonionic block has a uniform positive lower bound. For a unit vector \(\nu\), set
\begin{equation}\label{eq:display-033}
\Pi_\nu^{\rm oct}:=I_{16}-\cH(\nu\otimes\nu)_{\RR}.
\end{equation}
In a spin frame sending \(\nu\) to \((0,1)\), \eqref{eq:rank-one} and
\eqref{eq:equivariance} give
\begin{equation}\label{eq:display-034}
\Pi_\nu^{\rm oct}=\begin{pmatrix}I_8&0\\0&0\end{pmatrix},\qquad
\nu^\perp=\OO\oplus\operatorname{Im}\OO.
\end{equation}
Thus \(\Pi_\nu^{\rm oct}\) is the orthogonal projection onto the first
octonionic coordinate, whose eight real directions are tangent to the
boundary. The second octonionic coordinate contains the normal direction
and seven real tangent directions.

On the collar, we take \(\nu=Dd\).
For a real \(C^2\) function
\(\zeta\), define
\begin{equation}\label{eq:tangential-block}
\mathcal B_\nu[\zeta]:=\tr(\Pi_\nu^{\rm oct}D^2\zeta),
\qquad
\underline b:=\mathcal B_\nu[\uu],\qquad
\mathcal R_\nu:=-\mathcal B_\nu[d].
\end{equation}
The adapted frame gives \(\mathcal B_\nu[\zeta]=H(\zeta)_{11}\),
so \(\underline b\ge\kappa\) and \(\mathcal B_\nu[u]>0\).
Tracing the boundary identity \eqref{eq:tangential-identity} over
the range of \(\Pi_\nu^{\rm oct}\), we obtain
\begin{equation}\label{eq:block-boundary-identity}
\mathcal B_\nu[u]=\underline b-v_\nu\mathcal R_\nu
\qquad\text{on }\partial\Omega.
\end{equation}
On the boundary, \(\mathcal R_\nu\le0\) and \(v_\nu\ge0\) imply
\(\mathcal B_\nu[u]\ge\underline b\ge\kappa\).
The bounded derivatives of $d$ and $\underline u$ through order four
give uniform $C^2$ bounds for $\underline b$ and $\mathcal R_\nu$
by \eqref{eq:tangential-block}.

\begin{proposition}\label{prop:block-lower-bound}
There is a constant \(c_0>0\), depending only on the data, such that every
smooth admissible solution of \eqref{eq:path} satisfies
\begin{equation}\label{eq:block-lower-bound}
\mathcal B_\nu[u]\ge c_0\quad\text{on }\partial\Omega.
\end{equation}
\end{proposition}
\begin{proof}
Choose \(p\in\partial\Omega\) where the minimum
\begin{equation}\label{eq:display-inline-350}
\mu_0:=\min_{\partial\Omega}\mathcal B_\nu[u]>0
\end{equation}
is attained.
It suffices to consider \(\mu_0<\kappa/2\). At \(p\),
\eqref{eq:block-boundary-identity} gives
\begin{equation}\label{eq:display-035}
v_\nu(p)\mathcal R_\nu(p)=\underline b(p)-\mu_0>\frac{\kappa}{2}.
\end{equation}
By \eqref{eq:boundary-gradient}, we have
\begin{equation}\label{eq:display-inline-353}
\mathcal R_\nu(p)\ge\kappa/(2C_\nu).
\end{equation} The \(C^1\) bound for \(\mathcal R_\nu\) gives
a radius \(0<r<\delta\), depending only on the data, such that
\(\mathcal R_\nu\ge\kappa/(4C_\nu)\) in \(D:=\Omega\cap B_r(p)\). Set
\begin{equation}\label{eq:display-036}
\chi:=\frac{\mu_0-\underline b}{\mathcal R_\nu},\qquad
n:=\nu(p),\qquad W:=-\partial_n v-\chi.
\end{equation}
Since \(0<\mu_0<\kappa/2\), the \(C^2\) bounds for \(\underline b\) and
\(\mathcal R_\nu\) give \(\|\chi\|_{C^2(\overline D)}\le C\).
On \(\partial\Omega\cap\overline D\), the inequality \(v_\nu\ge0\) implies
\begin{equation}\label{eq:display-inline-366}
\partial_n v=(n\cdot\nu)v_\nu\le v_\nu.
\end{equation}
Thus
\begin{equation}\label{eq:display-037}
W\ge-v_\nu-\chi
=\frac{\mathcal B_\nu[u]-\mu_0}{\mathcal R_\nu}\ge0,
\qquad W(p)=0.
\end{equation}
Since \(n\) is fixed, \eqref{eq:first-differentiation} gives
\begin{equation}\label{eq:display-038}
\cL W=-n\cdot D_xg+\cL(\partial_n\uu)-\cL\chi.
\end{equation}
Hence \(|\cL W|\le C\mathcal T\), and \eqref{eq:first-order} gives
\(\|W\|_{C^0(\overline D)}\le C\). Applying
\eqref{eq:local-boundary-comparison} to \(W\), we obtain
\begin{equation}\label{eq:display-039}
-C\le W_\nu(p)=-v_{\nu\nu}(p)-\chi_\nu(p).
\end{equation}
The bounds for \(\chi_\nu\) and \(\underline u_{\nu\nu}\) then give us  
\begin{equation}\label{eq:normal-at-minimum}
u_{\nu\nu}(p)\le C.
\end{equation}

In a spin frame adapted to \(\nu(p)\), write
\begin{equation}\label{eq:normal-block-matrix}
H(u)=
\begin{pmatrix}
\mathcal B_\nu[u]&u_{\bar1 2}\\
\overline{u_{\bar1 2}}&u_{\bar2 2}
\end{pmatrix},
\qquad
u_{\bar2 2}
=u_{\nu\nu}+\sum_{j=1}^7u_{\eta_j\eta_j},
\end{equation}
where \(\eta_j:=(0,e_j)\) are unit tangent vectors.
Admissibility, \eqref{eq:tangential-bound}, and
\eqref{eq:normal-at-minimum} imply
\(0<u_{\bar2 2}(p)\le C_*\), with a fixed \(C_*>0\).
Evaluating the determinant equation at \(p\) gives
\begin{equation}\label{eq:minimum-closure}
\mu_0u_{\bar2 2}(p)
=g(t,p,u(p))^2+|u_{\bar1 2}(p)|^2\ge g_-^2,
\qquad \mu_0\ge\frac{g_-^2}{C_*}.
\end{equation}
Thus \eqref{eq:block-lower-bound} holds in both cases with
\begin{equation}\label{eq:display-inline-378}
c_0:=\min\{\kappa/2,g_-^2/C_*\}.
\end{equation}

\end{proof}

\subsection{The global second derivative estimate}
\begin{proposition}\label{prop:global-hessian}
Every smooth admissible solution of \eqref{eq:path} satisfies
\begin{equation}\label{eq:global-hessian}
\|u\|_{C^2(\overline\Omega)}\le C,\qquad
\lambda_-I_2\le H(u)\le\lambda_+I_2
\end{equation}
for fixed \(\lambda_-,\lambda_+>0\).
\end{proposition}
\begin{proof}
At each boundary point, choose a frame as in \eqref{eq:normal-block-matrix}.
Separating the normal coordinate \(x_{20}\) from the seven imaginary
coordinates gives
\begin{equation}\label{eq:offdiagonal-split}
u_{\bar1 2}
=\sum_{p=0}^7e_pu_{x_{1p}x_{20}}
+\sum_{p=0}^7\sum_{q=1}^7(e_p\bar e_q)u_{x_{1p}x_{2q}}.
\end{equation}
Proposition~\ref{prop:mixed} and \eqref{eq:tangential-bound}
bound the two sums, respectively. Hence \(|u_{\bar1 2}|\le C\). The determinant equation
and Proposition~\ref{prop:block-lower-bound} now give
\begin{equation}\label{eq:normal-block-upper}
0<u_{\bar2 2}
=\frac{g(t,q,u)^2+|u_{\bar1 2}|^2}{\mathcal B_\nu[u]}
\le\frac{g_+^2+C^2}{c_0}.
\end{equation}
By \eqref{eq:normal-block-matrix} and \eqref{eq:tangential-bound},
\begin{equation}\label{eq:display-040}
u_{\nu\nu}=u_{\bar2 2}-\sum_{j=1}^7u_{\eta_j\eta_j},
\qquad
-7C\le u_{\nu\nu}\le\frac{g_+^2+C^2}{c_0}+7C.
\end{equation}
Together with the tangential and mixed estimates, this gives
\begin{equation}\label{eq:boundary-hessian}
\sup_{\partial\Omega}|D^2u|\le C.
\end{equation}
Lemma~\ref{lem:second-derivatives} now gives \(D^2u\le CI_{16}\)
throughout \(\Omega\).
The eigenvalues \(\lambda_1^{\RR}\le\cdots\le\lambda_{16}^{\RR}\) of
\(D^2u\) satisfy
\begin{equation}\label{eq:real-eigenvalues}
\sum_{j=1}^{16}\lambda_j^{\RR}=\Delta u=\tr H(u)>0,\qquad
\lambda_1^{\RR}>-\sum_{j=2}^{16}\lambda_j^{\RR}\ge-15C.
\end{equation}
Thus \(|D^2u|\le C\).
By \eqref{eq:rank-one}, the inequality \(D^2u\le CI_{16}\) implies
\begin{equation}\label{eq:display-041}
0<H(u)\le\cH(CI_{16})=8CI_2.
\end{equation}
Take \(\lambda_+:=8C\). Since \(\det H(u)\ge g_-^2\),
\begin{equation}\label{eq:display-inline-390}
\lambda_{\min}(H(u))\ge g_-^2/\lambda_+.
\end{equation}

This proves \eqref{eq:global-hessian} with \(\lambda_-:=g_-^2/\lambda_+\).
\end{proof}

\section{Smooth solvability and translation quotients}\label{sec:existence}
\subsection{Higher-order estimates}
Estimates \eqref{eq:global-hessian} and \eqref{eq:ellipticity} give
\begin{equation}\label{eq:uniform-ellipticity}
\lambda_*I_{16}\le a\le\Lambda_*I_{16},\qquad
\lambda_*:=\frac{1}{4}\sqrt{\frac{\lambda_-}{\lambda_+}},\qquad
\Lambda_*:=\sqrt{\frac{\lambda_+}{\lambda_-}}.
\end{equation}
For the regularity estimates, extend \(\F\) to \(\Sym(16)\) by
\begin{equation}\label{eq:elliptic-extension}
\begin{aligned}
\cK&:=\left\{A\in\operatorname{Herm}_2(\OO):
\frac{\lambda_-}{2}I_2\le A\le2\lambda_+I_2\right\},\\
\widehat{\F}(S)&:=\inf_{A\in\cK}
\bigl\{\F_{\OO}(A)+D\F_{\OO}(A):(\cH(S)-A)\bigr\}.
\end{aligned}
\end{equation}
The extension is concave because it is an infimum of affine functions.
By \eqref{eq:support}, it agrees with \(\F\) whenever \(\cH(S)\in\cK\).
For \(A\in\cK\), the coefficient matrix of each affine function lies
between \(\lambda_*I_{16}\) and \(\Lambda_*I_{16}\) by
\eqref{eq:ellipticity}. Thus the extension satisfies
\begin{equation}\label{eq:extension-ellipticity}
\lambda_*\tr Q\le\widehat{\F}(S+Q)-\widehat{\F}(S)
\le\Lambda_*\tr Q\qquad(Q\ge0).
\end{equation}
By \eqref{eq:path-bounds} and \eqref{eq:global-hessian},
\begin{equation}\label{eq:composed-density}
\widehat{\F}(D^2u)=g(t,\cdot,u),\qquad
\|g(t,\cdot,u)\|_{C^{0,1/2}(\Om)}\le C.
\end{equation}
The operator \(\widehat{\F}\) is uniformly elliptic and concave.
Applying the Evans--Krylov interior estimate
\cite{Evans,KrylovInterior} and Krylov's boundary estimate
\cite{KrylovBoundary} to \eqref{eq:composed-density},
with \(u=\phi\) on \(\partial\Omega\), we obtain
\begin{equation}\label{eq:holder}
\|u\|_{C^{2,\alpha_*}(\Om)}\le C
\end{equation}
for some \(0<\alpha_*<1/2\), with \(\alpha_*\) and \(C\)
independent of \(t\).
See also \cite[Theorem~1.4(b)]{SilvestreSirakov}
for the global estimate. Then, by applying Schauder estimates to the original operator \(\F\), we get
\begin{equation}\label{eq:display-042}
\|u\|_{C^{k,\alpha_*}(\Om)}\le C_k\qquad(k\ge2),
\end{equation}
where \(C_k\) depends on the corresponding higher-order norms of the data.

\subsection{The continuity method and translation quotients}
\begin{proof}[Proof of Theorem~\ref{thm:main}]
Let
\begin{equation}\label{eq:zero-boundary-space}
\begin{aligned}
\mathcal X^{2,\alpha_*}
&:=\{w\in C^{2,\alpha_*}(\Om):w|_{\partial\Omega}=0\},\\
\mathcal A&:=\{w\in\mathcal X^{2,\alpha_*}:H(\uu+w)>0\text{ on }\Om\}.
\end{aligned}
\end{equation}
Then \(\mathcal A\) is open in \(\mathcal X^{2,\alpha_*}\). Define
\begin{equation}\label{eq:display-043}
\begin{aligned}
\mathcal M(t,w)&:=\F(D^2(\uu+w))-g(t,q,\uu+w),\\
\mathcal I&:=\{t\in[0,1]:\mathcal M(t,w)=0\text{ for some }w\in\mathcal A\}.
\end{aligned}
\end{equation}
Since \(\mathcal M(0,0)=0\), we have \(0\in\mathcal I\).
To prove openness, we linearize at a solution:
\begin{equation}\label{eq:display-044}
D_w\mathcal M(t,w)[\zeta]
=a^{ij}\zeta_{ij}-g_s(t,q,\uu+w)\zeta=\cL\zeta.
\end{equation}
The maximum principle, linear Dirichlet theory, and the Schauder
estimate show that
\begin{equation}\label{eq:isomorphism}
\cL:\mathcal X^{2,\alpha_*}\longrightarrow C^{\alpha_*}(\Om)
\quad\text{is an isomorphism and}\qquad
\|\zeta\|_{C^{2,\alpha_*}}\le C\|\cL\zeta\|_{C^{\alpha_*}}.
\end{equation}
The implicit function theorem shows that \(\mathcal I\) is open in \([0,1]\),
and the preceding regularity estimates give smooth solutions.

For closedness, take \(t_j\in\mathcal I\) with limit \(t_\infty\)
and corresponding solutions \(u_j\). By \eqref{eq:holder} and
\eqref{eq:global-hessian}, a subsequence satisfies
\begin{equation}\label{eq:display-045}
u_j\longrightarrow u_\infty\quad\text{in }C^2(\Om),\qquad
H(u_\infty)\ge\lambda_-I_2.
\end{equation}
The limit satisfies
\begin{equation}\label{eq:display-046}
\F(D^2u_\infty)=g(t_\infty,q,u_\infty),\qquad
u_\infty|_{\partial\Omega}=\phi.
\end{equation}
The uniform \(C^{2,\alpha_*}\) bound passes to the limit, which is
smooth by elliptic regularity. Hence \(\mathcal I=[0,1]\), and the uniqueness follows from Lemma~\ref{lem:comparison}.
\end{proof}

\begin{proof}[Proof of Corollary~\ref{cor:quotient}]
Let \(\pi:\OO^2\to X\) be the quotient map. Then we have
\begin{equation}\label{eq:display-047}
D^2(\zeta\circ\pi)=\pi^*D^2\zeta,\qquad
H(\zeta\circ\pi)=\pi^*H(\zeta).
\end{equation}
The sixteen parallel coordinate fields descend to \(X\).
Taking \(h\) to be the harmonic extension of \(\phi\) to \(M\) and setting
\begin{equation}\label{eq:display-inline-423}
p_0=\uu-\max_{\overline M}\uu-1,
\end{equation}
we obtain the comparison and derivative estimates on \(M\)
by the same argument as in Section~\ref{sec:estimates}.

For the boundary estimates, we cover \(\partial M\) by finitely
many coordinate neighborhoods in \(X\) and lift each to \(\RR^{16}\).
In each lifted chart, we choose a spin frame at a boundary point \(p\)
and use the affine field
\begin{equation}\label{eq:display-inline-429}
T(x)=\tau+K(x-p)
\end{equation}
for the mixed derivative estimate and the fixed direction \(\nu(p)\)
for the double normal derivative estimate.
Both estimates use the local barrier involving \(|x-p|^2\),
as in Section~\ref{sec:boundary}.
The finite cover makes the estimates uniform.

Since Dirichlet regularity and the isomorphism \eqref{eq:isomorphism}
also hold on \(M\), the same continuity argument gives existence.
Uniqueness follows from the comparison principle.
\end{proof}

\section{Subsolutions and annular domains}\label{sec:applications}
\subsection{Subsolutions from a defining function}
\begin{proposition}\label{prop:defining-function}
Suppose \(\Omega\) has a smooth defining function \(\rho\in C^\infty(\Om)\)
such that
\begin{equation}\label{eq:defining-function}
\rho<0\text{ in }\Omega,\qquad
\rho=0,\ d\rho\ne0\text{ on }\partial\Omega,\qquad
H(\rho)\ge\kappa_\rho I_2\text{ on }\Om
\end{equation}
for some \(\kappa_\rho>0\).
Then for each smooth boundary datum \(\phi\) and smooth density \(f\)
satisfying \eqref{eq:density}, there exists a subsolution as in
\eqref{eq:subsolution}.
\end{proposition}
\begin{proof}
Let \(\widetilde\phi\) be a smooth extension of \(\phi\) to \(\Om\).
For \(\Lambda_{\rm sub}\ge0\), set
\begin{equation}\label{eq:display-048}
\uu:=\widetilde\phi+\Lambda_{\rm sub}\rho,\qquad
C_\phi:=\sup_{\Om}\|H(\widetilde\phi)\|_{\rm op}.
\end{equation}
By \eqref{eq:defining-function},
\begin{equation}\label{eq:defining-function-bound}
H(\uu)\ge(\Lambda_{\rm sub}\kappa_\rho-C_\phi)I_2.
\end{equation}
Since \(\rho\le0\) and \(f_s\ge0\),
\begin{equation}\label{eq:defining-function-density}
\uu\le\widetilde\phi,\qquad
f(q,\uu)\le f(q,\widetilde\phi)\le
f_{\max}:=\max_{q\in\Om}f(q,\widetilde\phi(q)).
\end{equation}
Choose
\begin{equation}\label{eq:display-049}
\Lambda_{\rm sub}\ge
\frac{C_\phi+\sqrt{f_{\max}}+1}{\kappa_\rho}.
\end{equation}
Then
\begin{equation}\label{eq:defining-function-subsolution}
H(\uu)\ge(\sqrt{f_{\max}}+1)I_2,\qquad
\det H(\uu)>f_{\max}\ge f(q,\uu).
\end{equation}

The boundary equality \(\uu=\phi\) follows from \(\rho=0\) on \(\partial\Omega\).
\end{proof}

The bound \eqref{eq:defining-function-density} requires no growth
assumption on \(f(q,s)\) as \(s\to-\infty\).

\subsection{An annular family}
Let \(0<R_{\rm in}<R_{\rm out}\) and take
\begin{equation}\label{eq:annular-data}
\Omega:=\{q\in\OO^2:R_{\rm in}<|q|<R_{\rm out}\},\qquad
\phi:=\left.\frac{|q|^2}{16}\right|_{\partial\Omega},\qquad
f\equiv f_{\rm ann}\in(0,1].
\end{equation}
The function \(\uu:=|q|^2/16\) satisfies \(H(\uu)=I_2\) and
\begin{equation}\label{eq:display-inline-454}
\det H(\uu)=1\ge f_{\rm ann}.
\end{equation}

The unit normal pointing into the annulus is \(\nu=q/|q|\).
Near this sphere, \(d=|q|-R_{\rm in}\). Thus on \(|q|=R_{\rm in}\),
for any tangent vector \(\tau\),
\begin{equation}\label{eq:display-050}
D^2d(\tau,\tau)=\frac{|\tau|^2}{R_{\rm in}},\qquad
\mathcal R_\nu=-\frac{8}{R_{\rm in}}.
\end{equation}
For any defining function \(\rho\) negative in the annulus,
\(\rho_\nu<0\) on the inner sphere. Taking the trace of
\eqref{eq:tangential-identity} over the octonionic tangent plane gives
\begin{equation}\label{eq:annular-obstruction}
\mathcal B_\nu[\rho]=-\rho_\nu\mathcal R_\nu
=\frac{8\rho_\nu}{R_{\rm in}}<0
\qquad\text{on }\{|q|=R_{\rm in}\}.
\end{equation}
This excludes a strictly plurisubharmonic
defining function.

To solve the radial equation, write \(u(x)=\mathsf u(\vartheta)\), where
\begin{equation}\label{eq:display-051}
\vartheta:=\frac{|q|^2}{16},\qquad
\vartheta_{\rm in}:=\frac{R_{\rm in}^2}{16},\qquad
\vartheta_{\rm out}:=\frac{R_{\rm out}^2}{16}.
\end{equation}
In real coordinates \(x=q_{\RR}\), we have
\(\vartheta_i=x_i/8\) and \(\vartheta_{ij}=\delta_{ij}/8\). Hence
\begin{equation}\label{eq:radial-hessian}
u_{ij}=\frac{\mathsf u'}{8}\delta_{ij}+\frac{\mathsf u''}{64}x_ix_j,\qquad
H(u)=\mathsf u'I_2+\frac{\mathsf u''}{64}\cH(x\otimes x).
\end{equation}
By \eqref{eq:rank-one}, the eigenvalues of \(\cH(x\otimes x)\) are
\(0\) and \(|q|^2\). Since \(|q|^2/64=\vartheta/4\), the eigenvalues
of \(H(u)\) are
\begin{equation}\label{eq:radial-eigenvalues}
\mathsf u',\qquad \mathsf u'+\frac{\vartheta}{4}\mathsf u''.
\end{equation}
Their product gives the radial equation
\begin{equation}\label{eq:radial-equation}
\mathsf u'\left(\mathsf u'+\frac{\vartheta}{4}\mathsf u''\right)=f_{\rm ann}.
\end{equation}
Set \(y:=(\mathsf u')^2\). Since \(y'=2\mathsf u'\mathsf u''\),
\begin{equation}\label{eq:display-052}
y+\frac{\vartheta}{8}y'=f_{\rm ann},\qquad
(\vartheta^8y)'=8f_{\rm ann}\vartheta^7.
\end{equation}
Integration gives 
\begin{equation}\label{eq:display-inline-472}
y=f_{\rm ann}+C_{\rm rad}\vartheta^{-8}.
\end{equation}

The positive root and the inner boundary value determine
\begin{equation}\label{eq:radial-solution}
\begin{aligned}
\mathsf u'(\vartheta)&=\sqrt{f_{\rm ann}+C_{\rm rad}\vartheta^{-8}},\\
\mathsf u(\vartheta)&=\vartheta_{\rm in}
+\int_{\vartheta_{\rm in}}^\vartheta
\sqrt{f_{\rm ann}+C_{\rm rad}\theta^{-8}}\,d\theta.
\end{aligned}
\end{equation}
At \(\vartheta=\vartheta_{\rm out}\), the integral above depends continuously on \(C_{\rm rad}\ge0\)
and increases strictly from
\(\sqrt{f_{\rm ann}}(\vartheta_{\rm out}-\vartheta_{\rm in})\)
at \(C_{\rm rad}=0\) to infinity as \(C_{\rm rad}\to\infty\).
Thus, for \(0<f_{\rm ann}<1\), the outer boundary condition
determines a unique \(C_{\rm rad}>0\).
For \(f_{\rm ann}=1\), we take \(C_{\rm rad}=0\).
In both cases, we have
\begin{equation}\label{eq:display-053}
\mathsf u''=-\frac{4C_{\rm rad}\vartheta^{-9}}{\mathsf u'}\le0.
\end{equation}
By \eqref{eq:radial-eigenvalues} and \eqref{eq:radial-equation},
the eigenvalues in increasing order are
\begin{equation}\label{eq:display-054}
\lambda_1(H(u))=\frac{f_{\rm ann}}{\mathsf u'}>0,\qquad
\lambda_2(H(u))=\mathsf u'\ge\lambda_1(H(u)).
\end{equation}
Thus the radial solution is smooth and admissible on the closed annulus.
For \(f_{\rm ann}<1\), we have \(\mathsf u''<0\). Since
\(\mathsf u(\vartheta_{\rm in})=\vartheta_{\rm in}\) and
\(\mathsf u(\vartheta_{\rm out})=\vartheta_{\rm out}\), concavity
gives \(\mathsf u(\vartheta)\ge\vartheta\), as required by comparison.
We use the homogeneous radial equation in Section~\ref{sec:construction}
to determine the singularity of the Green function.

\section{The homogeneous equation and the punctured domains}\label{sec:construction}
We construct the Green function as a limit of solutions of
positive-density equations on punctured domains.

\subsection{The admissible branch and boundary geometry}
On the closure of the admissible cone of Proposition~\ref{prop:structure},
the homogeneous equation is
\begin{equation}\label{eq:mineigen}
\lambda_1(\cH(D^2u))=0 .
\end{equation}
Here $\lambda_1$ denotes the least Hermitian eigenvalue defined in
\eqref{eq:hermitian-eigenvalues}. For continuous $u$, we interpret
\eqref{eq:mineigen} in the viscosity sense: at every $x_0$, each real
$C^2$ function $\phi$ near $x_0$ with $\phi(x_0)=u(x_0)$ must satisfy
\begin{equation}\label{eq:display-055}
\begin{aligned}
\phi\ge u\text{ near }x_0
&\quad\Longrightarrow\quad H(\phi)(x_0)\ge0
&&\text{(upper test)},\\
\phi\le u\text{ near }x_0
&\quad\Longrightarrow\quad \lambda_1(H(\phi)(x_0))\le0
&&\text{(lower test)}.
\end{aligned}
\end{equation}
For a $C^2$ function $\phi$, the condition $H(\phi)(x_0)\ge0$ is
equivalent to $D^2\phi(x_0)$ having nonnegative trace on every octonionic
line through $x_0$ \cite[Proposition~3.1.8]{AleskerOct}. By the
restriction theorem \cite[Theorem~3.2]{HLrestriction}, an upper
semicontinuous function satisfies the upper test condition if and only
if its restriction to every affine octonionic line is subharmonic.

\begin{lemma}\label{lem:globalrho}
Suppose \(\Omega\Subset\OO^2\) has smooth strictly octonionic
pseudoconvex boundary. It admits a smooth defining function $\rho$
satisfying \eqref{eq:defining-function}. Moreover,
\begin{equation}\label{eq:Levimargin}
\mathcal R_\nu\ge r_0>0
\qquad\text{on }\partial\Omega.
\end{equation}
\end{lemma}
\begin{proof}
Let \(\rho_p\) be a local strictly plurisubharmonic defining function
near \(p\in\partial\Omega\), negative inside \(\Omega\). Its inward
normal derivative is negative, and the tangential trace in
\eqref{eq:tangential-identity} gives
\begin{equation}\label{eq:display-inline-504}
\mathcal B_\nu[\rho_p]=-(\rho_p)_\nu\mathcal R_\nu
\quad\text{on the boundary}.
\end{equation}
Thus \(\mathcal R_\nu>0\).
Compactness gives \eqref{eq:Levimargin}.

On a boundary collar, set \(\rho_0:=-d+Ad^2\). On the boundary,
\begin{equation}\label{eq:display-056}
H(\rho_0)=-H(d)+2A\cH(Dd\otimes Dd).
\end{equation}
At each boundary point, choose a spin frame with \(\nu=(0,1)\). Then
\begin{equation}\label{eq:display-057}
H(\rho_0)=
\begin{pmatrix}
\mathcal R_\nu&-H(d)_{12}\\
-H(d)_{21}&2A-H(d)_{22}
\end{pmatrix}.
\end{equation}
By \eqref{eq:Levimargin} and the smoothness of \(d\),
\begin{equation}\label{eq:display-058}
\mathcal R_\nu\ge r_0,\qquad |H(d)_{12}|+|H(d)_{22}|\le C,
\end{equation}
where \(C\) depends only on the boundary geometry. Choose \(A\) so that
\begin{equation}\label{eq:display-inline-511}
2A\ge C+(C^2+1)/r_0.
\end{equation}
Then
\begin{equation}\label{eq:rho0-determinant}
\begin{aligned}
\det H(\rho_0)
&=\mathcal R_\nu\bigl(2A-H(d)_{22}\bigr)-|H(d)_{12}|^2\\
&\ge r_0(2A-C)-C^2\ge1.
\end{aligned}
\end{equation}
Since \(H(\rho_0)_{11}\ge r_0\), this proves
\(H(\rho_0)>0\) on \(\partial\Omega\). By continuity and compactness,
we may shrink the collar so that
\begin{equation}\label{eq:display-059}
H(\rho_0)>0,\qquad 0\le d<\frac1{4A}.
\end{equation}
In this collar,
\begin{equation}\label{eq:display-060}
\begin{gathered}
1-2Ad\ge\frac12,\qquad D\rho_0=-(1-2Ad)Dd\ne0,\\
\rho_0=-d(1-Ad)<0\qquad(d>0).
\end{gathered}
\end{equation}

To extend the defining function to the interior, choose \(c>0\)
so that the levels \(\rho_0=-2c,-c/2\) lie in the collar.
Since \(\Omega\) is bounded, we can choose a quadratic function
\begin{equation}\label{eq:interior-quadratic}
Q(x)=\eta|x|^2-C_0,\qquad \eta>0,\qquad
-3c/2<Q<-c\quad\text{on }\Om.
\end{equation}

We join \(\rho_0\) to \(Q\) between these levels by the regularized
maximum
\begin{equation}\label{eq:display-061}
M_\sigma(s,t):=\tfrac12\bigl(s+t+\chi_\sigma(s-t)\bigr),
\end{equation}
where \(\chi_\sigma\) is smooth and convex,
\(|\chi_\sigma'|\le1\), and \(\chi_\sigma(z)=|z|\) for \(|z|\ge\sigma\).
For real functions \(v,w\), \eqref{eq:product-chain} gives
\begin{equation}\label{eq:regmax}
\begin{aligned}
H(M_\sigma(v,w))
&=\theta H(v)+(1-\theta)H(w)\\
&\quad+\tfrac12\chi_\sigma''(v-w)
\cH(D(v-w)\otimes D(v-w)),\qquad
\theta:=\tfrac12(1+\chi_\sigma'(v-w)).
\end{aligned}
\end{equation}
Since \(0\le\theta\le1\) and the last term in \eqref{eq:regmax}
is nonnegative, the regularized maximum preserves strict
plurisubharmonicity. Choose
\(0<\sigma<c/2\). At the two ends of the transition region,
\begin{equation}\label{eq:display-062}
\begin{aligned}
\rho_0=-c/2&\quad\Longrightarrow\quad \rho_0-Q>c/2>\sigma,\\
\rho_0=-2c&\quad\Longrightarrow\quad Q-\rho_0>c/2>\sigma.
\end{aligned}
\end{equation}
Thus \(M_\sigma(\rho_0,Q)\) agrees with \(\rho_0\) near the outer
level and with \(Q\) near the inner one. Since
\begin{equation}\label{eq:display-inline-534}
0\le M_\sigma(s,t)-\max\{s,t\}\le\sigma/2,
\end{equation}
we also have
\begin{equation}\label{eq:display-063}
M_\sigma(\rho_0,Q)
\le\max\{\rho_0,Q\}+\frac\sigma2
\le-\frac c2+\frac\sigma2<0
\qquad(-2c\le\rho_0\le-c/2).
\end{equation}
Extend this function by \(\rho_0\) near the boundary and by \(Q\)
in the interior. The resulting function is smooth, negative in
\(\Omega\), and strictly plurisubharmonic by \eqref{eq:regmax}.
Compactness gives the uniform lower bound in \eqref{eq:defining-function}.
\end{proof}

\subsection{The pole and a strict subsolution}
After a translation, we may place the pole at $0$. By
\eqref{eq:radial-eigenvalues}, the homogeneous radial equation with
$\mathsf u'>0$ reduces to $\vartheta\mathsf u''+4\mathsf u'=0$.
Its singular part is proportional to $\vartheta^{-3}$. Normalizing
the pole coefficient to one and using \eqref{eq:radial-hessian}, we obtain
\begin{equation}\label{eq:pole}
\Phi(x):=-|x|^{-6},\qquad
H(\Phi)=48|x|^{-8}I_2-48|x|^{-10}\cH(x\otimes x).
\end{equation}
By \eqref{eq:rank-one}, its eigenvalues are $0$ and $48|x|^{-8}$.
Thus $\Phi$ solves \eqref{eq:mineigen} away from $0$. See also
\cite[Proposition~5.1]{Wang}.

\begin{proposition}\label{prop:singularsub}
The defining function $\rho$ can be scaled so that $H(\rho)\ge I_2$
and there is a smooth subsolution $\uu$ on $\Om\setminus\{0\}$ such that
\begin{equation}\label{eq:singularsub}
\begin{gathered}
\uu=\Phi+|x|^2+B\quad(0<|x|\le R_1),\qquad
\uu=\rho-1\quad\text{near }\partial\Omega,\\
\uu\le-1,\qquad H(\uu)\ge I_2,\qquad |\uu-\Phi|\le C
\end{gathered}
\end{equation}
for some $B\in\RR$ and $R_1>0$.
\end{proposition}
\begin{proof}
Choose $B$ so that 
\begin{equation}\label{eq:singular-quadratic-bound}
w:=\Phi+|x|^2+B\le-3
\quad\text{throughout }\Om\setminus\{0\}.
\end{equation}
Then 
\begin{equation}\label{eq:display-inline-555}
H(w)=H(\Phi)+16I_2\ge16I_2.
\end{equation}

Choose $R_1>0$ with $\overline B_{2R_1}\subset\Omega$. By \eqref{eq:defining-function},
scale $\rho$ so that $H(\rho)\ge I_2$ on
$\Om$ and $\rho-1<w-1$ on $R_1\le|x|\le2R_1$.
Near $\partial\Omega$ we have $\rho-1>w+1$.

For $0<\sigma<1$, set
\begin{equation}\label{eq:display-064}
\uu(x):=
\begin{cases}
w(x),&0<|x|\le R_1,\\
M_\sigma(w(x),\rho(x)-1),&|x|>R_1.
\end{cases}
\end{equation}
The two formulas agree near $|x|=R_1$, and $\uu=\rho-1$ near
$\partial\Omega$. Where the regularized maximum is used,
\eqref{eq:regmax} gives
\begin{equation}\label{eq:display-065}
H(\uu)\ge\theta H(w)+(1-\theta)H(\rho)
\ge(1+15\theta)I_2\ge I_2.
\end{equation}
The same lower bound holds in the inner ball. Where
$|w-(\rho-1)|<\sigma$,
\begin{equation}\label{eq:display-066}
\max\{w,\rho-1\}\le-2,\qquad
\uu\le-2+\frac\sigma2<-1.
\end{equation}
Elsewhere, the regularized maximum equals
\begin{equation}\label{eq:display-inline-570}
\max\{w,\rho-1\}\le-1.
\end{equation}
Finally, the difference satisfies
\begin{equation}\label{eq:display-inline-571}
\uu-\Phi=|x|^2+B
\end{equation}
near $0$ and is bounded on the remaining compact set.
\end{proof}

For the approximating problems, fix
\begin{equation}\label{eq:display-inline-573}
0<r_*<\min\{1,R_1/4\}
\end{equation}
and take
\begin{equation}\label{eq:rectangle}
0<r<r_*,\qquad 0<\varepsilon\le1,\qquad
\Omega_r:=\Omega\setminus\overline B_r,\qquad c_r:=-r^{-6}+r^2+B .
\end{equation}
The subsolution in \eqref{eq:singularsub} takes the values $-1$ and
$c_r$ on $\partial\Omega$ and $\partial B_r$, respectively, and satisfies
$\det H(\uu)\ge1$.
Theorem~\ref{thm:main} therefore gives a smooth admissible solution
$u=u^{\varepsilon,r}$ of
\begin{equation}\label{eq:approxproblem}
\det H(u)=\varepsilon\quad\text{in }\Omega_r,\qquad
u=-1\quad\text{on }\partial\Omega,\qquad
u=c_r\quad\text{on }\partial B_r .
\end{equation}
For each fixed pair $(\varepsilon,r)$, the linearization is elliptic
on $\overline\Omega_r$.

\section{Bounds for the solution and its gradient on punctured domains}\label{sec:greenC1}
\subsection{Pointwise bounds and boundary gradient estimates}
All constants below are independent of \(\varepsilon\) and \(r\) in
\eqref{eq:rectangle}. Comparison with \(\uu\) and the maximum principle
for \(\Delta u>0\) give
\begin{equation}\label{eq:display-067}
\uu\le u\le-1.
\end{equation}
Choose \(C_+\) so that
\(\Phi+C_+\ge-1\) on \(\partial\Omega\) and
\(C_+\ge r_*^2+B\). Then \(\Phi+C_+\ge u\) on both boundary components.
If \(u-(\Phi+C_+)\) had a positive interior maximum,
\eqref{eq:rank-one} would give
\begin{equation}\label{eq:display-068}
D^2u\le D^2\Phi\quad\Longrightarrow\quad 0<H(u)\le H(\Phi),
\end{equation}
which contradicts \(\lambda_1(H(\Phi))=0\). Thus
\(u\le\Phi+C_+\). Together with \eqref{eq:singularsub}, this gives
\begin{equation}\label{eq:greenheight}
|u-\Phi|\le C,\qquad c|x|^{-6}\le -u(x)\le C|x|^{-6}.
\end{equation}
For the lower bound on \(-u\), use the first estimate in
\eqref{eq:greenheight} near the pole and \(-u\ge1\) elsewhere.
The upper bound follows from the same estimate because \(\Omega\) is bounded.

For the gradient bound on the outer boundary, decrease $R_1$,
then $r_*$, so that $\Phi+C_+\le-2$ on
$\partial B_{R_1}$ and $r_*<R_1/4$. Let $h$ be harmonic in
$\Omega\setminus\overline B_{R_1}$, with values $-1$ on
$\partial\Omega$ and $-2$ on $\partial B_{R_1}$. Comparison gives
$u\le h$. Since the fixed domain is connected, the strong maximum principle
gives $h<-1$ in the interior.
On \(\partial\Omega\), \(\uu=u=h=-1\). Comparing inward normal
derivatives and applying the Hopf lemma uniformly on the compact boundary gives
\begin{equation}\label{eq:outergradient}
-C\le (\uu)_\nu\le u_\nu\le h_\nu\le-c<0
\qquad\text{on }\partial\Omega.
\end{equation}

On the inner boundary, $u=\uu=c_r$. For the upper barrier, take
\begin{equation}\label{eq:display-069}
h_r(x):=c_r+2(r^{-6}-|x|^{-6}),\qquad H(h_r)=2H(\Phi).
\end{equation}
Then \(h_r=c_r\) on \(\partial B_r\) and, for \(r_*\) small,
\begin{equation}\label{eq:display-inline-617}
h_r=r^{-6}+r^2+B-2|x|^{-6}\ge-1
\quad\text{on }\partial\Omega.
\end{equation}
Since \(H(u)>0\) and \(\lambda_1(H(h_r))=0\), the same comparison
argument gives \(u\le h_r\). We also have $\uu\le u$ by comparison.
Differentiating at \(|x|=r\) along \(\nu=x/r\) gives
\begin{equation}\label{eq:innergradient}
6r^{-7}+2r=\uu_\nu\le u_\nu\le (h_r)_\nu=12r^{-7}.
\end{equation}
The tangential gradient vanishes on $\partial\Omega_r$.
Thus \eqref{eq:greenheight}, \eqref{eq:outergradient}, and
\eqref{eq:innergradient} give
\begin{equation}\label{eq:weightedboundaryC1}
|x|\,|Du(x)|\le C(-u(x))\qquad(x\in\partial\Omega_r).
\end{equation}

\subsection{Spin rotations and dilations}
To extend the gradient estimate to the interior, we differentiate along
spin rotations and dilations. Given unit vectors $e,\xi$, take $\beta=e\cdot\xi$. By
\eqref{eq:quantitative-transitivity} we can choose $K\in\spin(9)$
with $Ke=\xi-\beta e$ and $|K|\le C$. Hence
\begin{equation}\label{eq:generator}
M:=\beta I_{16}+K,\qquad Me=\xi,\qquad |M|\le C .
\end{equation}
For each such $M$, set
\begin{equation}\label{eq:flow}
u_t(x):=e^{-2\beta t}u(e^{tM}x).
\end{equation}
At each interior point this is defined for small $|t|$. Since
$e^{tM}=e^{\beta t}e^{tK}$ and $e^{tK}\in\Spin(9)$, differentiation
in $x$ and \eqref{eq:equivariance} give
\begin{equation}\label{eq:flowequation}
D^2u_t(x)=e^{-tK}D^2u(e^{tM}x)e^{tK},\qquad
\F(D^2u_t)=\sqrt\varepsilon .
\end{equation}
With $L$ as in \eqref{eq:linearized-operator}, we have
$Lu=\sqrt\varepsilon$ by \eqref{eq:mixed-linearization}.
Differentiating \eqref{eq:flowequation} at $t=0$ gives
\begin{equation}\label{eq:firstflow}
V_M:=\left.\frac{\partial u_t}{\partial t}\right|_{t=0}
=Du\cdot Mx-2\beta u,\qquad LV_M=0.
\end{equation}
The expression in
\eqref{eq:firstflow} defines \(V_M\) throughout \(\overline\Omega_r\),
even when \(e^{tM}\) does not preserve \(\Omega_r\).
By \eqref{eq:weightedboundaryC1} and \eqref{eq:generator}, we can
choose \(A_0>0\), uniformly in \(M\), such that
\begin{equation}\label{eq:display-070}
\begin{aligned}
\pm V_M+A_0u&\le0&&\text{on }\partial\Omega_r,\\
L(\pm V_M+A_0u)&=A_0\sqrt\varepsilon>0&&\text{in }\Omega_r.
\end{aligned}
\end{equation}
The maximum principle gives \(|V_M|\le A_0(-u)\) in \(\Omega_r\).

Fix \(x\in\Omega_r\) and a real unit vector \(\xi\). Choose \(M\) as in \eqref{eq:generator} with
\(e=x/|x|\). Then 
\begin{equation}\label{eq:display-inline-651}
Du\cdot Mx=|x|\partial_\xi u,
\end{equation}
 so
\eqref{eq:firstflow} gives
\begin{equation}\label{eq:display-071}
|x|\,|\partial_\xi u|
=|V_M+2\beta u|\le(A_0+2)(-u).
\end{equation}
Since the constants are uniform in \(M\), \eqref{eq:greenheight} yields
\begin{equation}\label{eq:weightedC1}
|x|\,|Du(x)|\le C(-u(x)),\qquad |Du(x)|\le C|x|^{-7}.
\end{equation}

\section{Boundary Hessian estimates for the punctured problem}\label{sec:uniformboundary}
We use the boundary argument of Section~\ref{sec:boundary}, keeping
the constants independent of a positive lower bound for the density.

Fix a point $p$ on a smooth boundary patch. Suppose that $u$ has constant
boundary value on the patch, $\F(D^2u)$ is a positive constant bounded
by one, and $u$ has a bounded $C^1$ norm in the collar. Let $b\le u$
there, with $b=u$ on the patch, $H(b)\ge\kappa I_2$ and
$\|b\|_{C^2}\le C$.
Lemma~\ref{lem:collar}, applied to $u-b$, gives
\eqref{eq:local-boundary-comparison} with constants depending only on
these bounds and the fixed boundary charts. Indeed,
\eqref{eq:collar-determinant}--\eqref{eq:collar-linearization}
use only the upper bound for $\F(D^2u)$.

For the affine field $T$ of Proposition~\ref{prop:mixed}, the constant
density and \eqref{eq:invariance-differentiated}, together with the
constant boundary value, give
\begin{equation}\label{eq:display-072}
L(Tu)=0\quad\text{in the collar},\qquad
Tu=u_\nu T\cdot\nu=O(|x-p|^2)\quad\text{on the boundary patch}.
\end{equation}
Applying \eqref{eq:local-boundary-comparison} to $\pm Tu$, we obtain
\begin{equation}\label{eq:display-073}
-C\le\partial_\nu(\pm Tu)(p)=\pm u_{\nu\tau}(p),
\end{equation}
where the equality follows from 
\begin{equation}\label{eq:display-inline-666}
(K\nu)\cdot Du(p)=0.
\end{equation}

Thus $|u_{\nu\tau}(p)|\le C$, and
\eqref{eq:tangential-identity} bounds the tangential second derivatives.
To estimate the double normal derivative, suppose
$\mathcal B_\nu[u]\ge c>0$ on the patch. Then
\eqref{eq:offdiagonal-split} gives $|u_{\bar1 2}|\le C$, and the
determinant equation yields
\begin{equation}\label{eq:display-074}
0<u_{\bar2 2}
=\frac{\F(D^2u)^2+|u_{\bar1 2}|^2}{\mathcal B_\nu[u]}
\le\frac{1+C^2}{c}.
\end{equation}
By \eqref{eq:normal-block-matrix},
\begin{equation}\label{eq:display-075}
u_{\nu\nu}=u_{\bar2 2}-\sum_{j=1}^7u_{\eta_j\eta_j},\qquad
|D^2u|\le C
\end{equation}
on the patch. The constants in this estimate are independent of any
positive lower bound for $\F(D^2u)$.

Apply this estimate in a fixed collar of $\partial\Omega$,
with $b=\uu=\rho-1$.
The required $C^1$ bound follows from \eqref{eq:greenheight} and
\eqref{eq:weightedC1}. Since $u=-1$ on $\partial\Omega$,
\eqref{eq:tangential-identity}, \eqref{eq:Levimargin} and
\eqref{eq:outergradient} give
\begin{equation}\label{eq:outerC2}
\mathcal B_\nu[u]=-u_\nu\mathcal R_\nu\ge cr_0,
\qquad |D^2u|\le C\quad\text{on }\partial\Omega .
\end{equation}

For the inner boundary, we rescale to a fixed annulus. Set
\begin{equation}\label{eq:display-inline-676}
U(y):=r^6(u(ry)-c_r),\qquad 1\le |y|\le2.
\end{equation}
Then $U=0$ on $|y|=1$, and
\begin{equation}\label{eq:display-076}
D_yU=r^7D_xu(ry),\qquad D_y^2U=r^8D_x^2u(ry),\qquad
\det H(U)=r^{16}\varepsilon .
\end{equation}
Equations \eqref{eq:greenheight} and \eqref{eq:weightedC1} give
$\|U\|_{C^1}\le C$ on this fixed annulus.
Since $2r<R_1$, comparison with \eqref{eq:singularsub} gives
\begin{equation}\label{eq:display-077}
U(y)\ge 1-|y|^{-6}+r^8(|y|^2-1)
\ge \frac{|y|^2-1}{4}=:b(y),\qquad H(b)=4I_2 .
\end{equation}
Here 
\begin{equation}\label{eq:display-inline-682}
1-|y|^{-6}=(|y|^2-1)(|y|^{-2}+|y|^{-4}+|y|^{-6}).
\end{equation}

On $|y|=1$, the inward normal is $\nu=y$. By
\eqref{eq:tangential-identity} and \eqref{eq:innergradient}, we have
\begin{equation}\label{eq:display-078}
\mathcal B_\nu[U]=8U_\nu=8r^7u_\nu(ry)\ge48 .
\end{equation}
Applying \eqref{eq:display-075}
to \(U\), we obtain \(|D_y^2U|\le C\) on \(|y|=1\). Then, 
returning to the original variables, we get
\begin{equation}\label{eq:innerC2}
|D^2u(x)|\le Cr^{-8}\qquad (|x|=r).
\end{equation}
Together with \eqref{eq:outerC2} and \eqref{eq:greenheight}, this gives us 
\begin{equation}\label{eq:weightedboundaryC2}
|x|^2|D^2u(x)|\le C(-u(x))\qquad(x\in\partial\Omega_r),
\end{equation}
with $C$ independent of $r$ and $\varepsilon$.

\section{The weighted Hessian estimate and the Green function}\label{sec:limit}
\subsection{Second variation}
For the flow \eqref{eq:flow}, we define
\begin{equation}\label{eq:display-inline-690}
W_M:=\left.\partial_t^2u_t\right|_{t=0}.
\end{equation}
Differentiating the flow twice, we get 
\begin{equation}\label{eq:secondvariation}
W_M=D^2u(Mx,Mx)+Du\cdot M^2x-4\beta Du\cdot Mx+4\beta^2u .
\end{equation}
To obtain a differential inequality for \(W_M\), we use
\eqref{eq:flowequation}, which holds near each interior point
for small \(|t|\). Differentiating
\(\F(D^2u_t)=\sqrt\varepsilon\) twice at \(t=0\)
and using the concavity of \(\F\), we obtain
\begin{equation}\label{eq:secondflowPDE}
LW_M=-D^2\F(D^2u)[D^2V_M,D^2V_M]\ge0,
\end{equation}
where $V_M$ is the first variation in \eqref{eq:firstflow}.

On $\partial\Omega_r$, \eqref{eq:generator}, \eqref{eq:weightedC1} and
\eqref{eq:weightedboundaryC2} bound $|W_M|$ by $C(-u)$.
Choose $A_2$ large enough, independently of $M$, so that
\begin{equation}\label{eq:display-079}
\begin{aligned}
W_M+A_2u&\le0&&\text{on }\partial\Omega_r,\\
L(W_M+A_2u)&\ge A_2\sqrt\varepsilon>0&&\text{in }\Omega_r.
\end{aligned}
\end{equation}
The maximum principle gives
\begin{equation}\label{eq:display-080}
W_M\le A_2(-u)\qquad\text{in }\Omega_r.
\end{equation}
For each unit $\xi$, choose $M$ with $M(x/|x|)=\xi$ as in
\eqref{eq:generator}. Using \eqref{eq:weightedC1} in \eqref{eq:secondvariation} gives
\begin{equation}\label{eq:display-081}
\begin{aligned}
|x|^2D^2u(\xi,\xi)
&=W_M-Du\cdot M^2x+4\beta Du\cdot Mx-4\beta^2u\\
&\le C(-u).
\end{aligned}
\end{equation}
Since $\Delta u>0$, the argument in \eqref{eq:real-eigenvalues}
also bounds the real eigenvalues from below. Hence, by
\eqref{eq:greenheight},
\begin{equation}\label{eq:fullHessian}
|D^2u^{\varepsilon,r}(x)|\le C(-u^{\varepsilon,r}(x))|x|^{-2}
\le C|x|^{-8}
\end{equation}
for all parameters in \eqref{eq:rectangle}.

\subsection{Passage to the homogeneous equation}
Let $(\varepsilon_j,r_j)\to(0,0)$ in \eqref{eq:rectangle}
and set $u_j:=u^{\varepsilon_j,r_j}$. On each compact subset of
$\Om\setminus\{0\}$, the functions $u_j$ are defined for large
$j$, and estimates \eqref{eq:greenheight}, \eqref{eq:weightedC1} and
\eqref{eq:fullHessian} bound $u_j$ and its first two derivatives uniformly. In a fixed boundary chart, the second derivatives
of the pullbacks are linear combinations of $D^2u_j$ and $Du_j$
with bounded coefficients, so the bounds also hold on half-balls
at the outer boundary. By compactness and a diagonal argument,
we can pass to a subsequence for which
\begin{equation}\label{eq:display-082}
\begin{aligned}
u_j+1&\longrightarrow G
&&\text{in }C^{1,\theta}_{\rm loc}(\Om\setminus\{0\}),\quad 0<\theta<1,\\
D^2u_j&\overset{*}{\rightharpoonup}D^2G
&&\text{in }L^\infty_{\rm loc}(\Omega\setminus\{0\}).
\end{aligned}
\end{equation}
The Hessian bound passes to the limit. Hence
$G\in C^{1,1}_{\rm loc}(\Om\setminus\{0\})$, and
\eqref{eq:weightedC1} and \eqref{eq:fullHessian} give \eqref{eq:mainreg}.
By \eqref{eq:approxproblem} and \eqref{eq:greenheight},
\begin{equation}\label{eq:display-083}
G|_{\partial\Omega}=0,\qquad G(x)=-|x|^{-6}+O(1)\quad(x\to0).
\end{equation}
Set $G(0)=-\infty$.

Since $H(u_j)>0$ and $\det H(u_j)=\varepsilon_j$,
\begin{equation}\label{eq:display-inline-715}
0<\lambda_1(H(u_j))\le\sqrt{\varepsilon_j}.
\end{equation}
The map $S\mapsto\lambda_1(\cH(S))$ is continuous and nondecreasing
under addition of a positive semidefinite matrix. Stability of viscosity
solutions on relatively compact subdomains of $\Omega\setminus\{0\}$
therefore gives $\lambda_1(H(G))=0$ in the sense of \eqref{eq:mineigen}.
Thus $G$ is octonionic plurisubharmonic there. On each octonionic
line through $0$, its restriction is subharmonic off $0$ and bounded
above near $0$. The removable singularity theorem extends this
restriction subharmonically across $0$, with value
\begin{equation}\label{eq:display-inline-724}
\limsup_{x\to0}G(x)=-\infty.
\end{equation}
Hence $G$ is octonionic plurisubharmonic on $\Omega$.

To identify the Monge--Ampère measure, take standard mollifications $v_j$ of a locally
$C^{1,1}$ octonionic plurisubharmonic function $v$. On each compact subset $K$ of its domain, for all large $j$,
\begin{equation}\label{eq:display-084}
\begin{gathered}
v_j\longrightarrow v\quad\text{uniformly on }K,\qquad
D^2v_j\longrightarrow D^2v\quad\text{a.e.},\\
\|D^2v_j\|_{L^\infty(K)}\le C_K,\qquad H(v_j)\ge0.
\end{gathered}
\end{equation}
Dominated convergence gives
\begin{equation}\label{eq:mollified-density}
\det H(v_j)\longrightarrow\det H(v)
\quad\text{in }L^1_{\rm loc},\qquad H(v)\ge0\quad\text{a.e.}
\end{equation}
Weak continuity \cite[Theorem~3.2.7]{AleskerOct} then identifies the
Monge--Ampère measure of $v$ with $\det H(v)\,dx$.

Applying this density formula to the limit and using weak continuity
along $u_j+1\to G$ on compact subsets away from the pole, we obtain
\begin{equation}\label{eq:display-085}
\det H(G)\,dx=\lim_{j\to\infty}\varepsilon_j\,dx=0
\qquad\text{weakly in }\Omega\setminus\{0\}.
\end{equation}

\subsection{Uniqueness}
Suppose $G_1,G_2$ satisfy \eqref{eq:greenproblem} with pole $0$.
For each $\eta>0$, the positive multiple $(1+\eta)G_1$ is
a solution of \eqref{eq:mineigen} off $0$ by homogeneity.
As $x\to0$,
\begin{equation}\label{eq:display-086}
(1+\eta)G_1(x)-G_2(x)=-\eta|x|^{-6}+O(1)<0 .
\end{equation}
On the outer boundary both functions vanish. The comparison principle
for the closed admissible cone \cite{HLdual,HLsing}
on $\Omega\setminus\overline B_s$ gives $(1+\eta)G_1\le G_2$ there
for all sufficiently small $s$. Letting $s\downarrow0$, then
$\eta\downarrow0$, gives $G_1\le G_2$. Interchanging $G_1$ and
$G_2$ proves uniqueness.

Equations \eqref{eq:display-067} and \eqref{eq:display-082} give $G\le0$.
The asymptotic formula \eqref{eq:display-083} rules out the constant zero
function, so the strong maximum principle for real subharmonic functions
gives $G<0$ off $0$.
Local compactness and uniqueness give convergence of the full family as
$(\varepsilon,r)\to(0,0)$: $u^{\varepsilon,r}+1\to G$ in
$C^{1,\theta}$ on each compact subset
of $\Om\setminus\{0\}$ for $0<\theta<1$. Translation gives the bounds
of Theorem~\ref{thm:green} for any pole $a\in\Omega$.

\section{Approximation error and asymptotic profile}\label{sec:profile}
\subsection{The approximation error}
\begin{proposition}\label{prop:rate}
Let $a=0$ and take $\rho$ as in Proposition~\ref{prop:singularsub}.
For some $C>0$ and all sufficiently small $r$, the solutions of
\eqref{eq:approxproblem} satisfy
\begin{equation}\label{eq:orderingrate}
(1+Cr^6)G+\sqrt\varepsilon\,\rho
\le u^{\varepsilon,r}+1\le(1-Cr^6)G
\quad\text{on }\overline{\Omega_r}.
\end{equation}
In particular, for each compact $K\subset\Om\setminus\{0\}$,
\begin{equation}\label{eq:rate}
\|u^{\varepsilon,r}+1-G\|_{C^0(K)}\le C_K(r^6+\sqrt\varepsilon)
\qquad (r<\min\{r_*,\operatorname{dist}(0,K)\}).
\end{equation}
\end{proposition}
\begin{proof}
On $\partial B_r$, \eqref{eq:greenproblem} and \eqref{eq:approxproblem}
give 
\begin{equation}\label{eq:display-inline-766}
|u^{\varepsilon,r}+1-G|\le C_0
\end{equation}
and $-G\ge r^{-6}/2$.
Choose $C\ge2C_0$ and decrease $r_*$ so that $Cr_*^6<1/2$.
On the inner sphere,
\begin{equation}\label{eq:display-087}
Cr^6(-G)\ge C/2\ge C_0\ge|u^{\varepsilon,r}+1-G|.
\end{equation}
Together with $\rho\le0$, this gives \eqref{eq:orderingrate} there.
All three functions vanish on the outer boundary.

The upper comparison function $(1-Cr^6)G$ has zero Monge--Ampère
measure. For the lower comparison function 
\begin{equation}\label{eq:display-inline-773}
(1+Cr^6)G+\sqrt\varepsilon\,\rho,
\end{equation}
we use $H(G)\ge0$ almost everywhere and $H(\rho)\ge I_2$ from
Proposition~\ref{prop:singularsub}. Thus
\begin{equation}\label{eq:display-088}
H\big((1+Cr^6)G+\sqrt\varepsilon\,\rho\big)
\ge \sqrt\varepsilon I_2\quad\text{almost everywhere in }\Omega_r.
\end{equation}
Both comparison functions are continuous and octonionic
plurisubharmonic on the punctured domain. By the density formula in
Section~\ref{sec:limit}, the lower comparison function has measure at least
$\varepsilon\,dx$, the measure of $u^{\varepsilon,r}+1$.
The comparison principle
\cite[Proposition~4.3]{Wang} gives \eqref{eq:orderingrate}.
Since $G\le0$ and $\rho\le0$, \eqref{eq:orderingrate} gives
\begin{equation}\label{eq:display-089}
|u^{\varepsilon,r}+1-G|\le Cr^6(-G)+\sqrt\varepsilon(-\rho).
\end{equation}
Taking the supremum on $K$ proves \eqref{eq:rate}.
\end{proof}

\subsection{Sharpness on the ball}
On the ball \(B_R(a)\), \eqref{eq:pole} and uniqueness give
\begin{equation}\label{eq:ballgreen}
G_{B_R(a)}(q,a)=R^{-6}-|q-a|^{-6}.
\end{equation}
With \(x=q-a\) and \(r=|x|\), its real derivatives are
\begin{equation}\label{eq:display-090}
DG=6r^{-8}x,\qquad D^2G=6r^{-8}I_{16}-48r^{-10}x\otimes x.
\end{equation}
The radial eigenvalue is \(-42r^{-8}\) and the fifteen tangential
eigenvalues are \(6r^{-8}\). Thus \(|DG|=6r^{-7}\) and
\begin{equation}\label{eq:display-inline-787}
\|D^2G\|_{\mathrm{op}}=42r^{-8},
\end{equation}
so both exponents in
\eqref{eq:mainreg} are sharp.

\subsection{The rescaled profile}
\begin{proof}[Proof of Corollary~\ref{cor:profile}]
Fix a compact \(K\subset\RR^{16}\setminus\{0\}\) and choose
\(d>0\) so that the closed \(2d\)-neighborhood \(K'\) of \(K\)
is disjoint from \(0\).
For sufficiently small \(\tau>0\), \(a+\tau K'\subset\Omega\).
Set 
\begin{equation}\label{eq:display-inline-796}
w_\tau(y):=G_\tau(y)+|y|^{-6}.
\end{equation}
By \eqref{eq:greenproblem}
and \eqref{eq:mainreg},
\begin{equation}\label{eq:display-091}
\|w_\tau\|_{C^0(K')}\le C\tau^6,\qquad
\|D^2w_\tau\|_{L^\infty(K')}\le C_{K'}.
\end{equation}
The second bound follows from
\begin{equation}\label{eq:display-inline-797}
D_y^2G_\tau(y)=\tau^8D_q^2G_\Omega(a+\tau y,a)
\quad\text{almost everywhere}
\end{equation}
and \eqref{eq:mainreg}, which gives the bound \(C|y|^{-8}\)
after rescaling.
For \(y\in K\), a unit vector \(e\), and \(0<h<d\), the
first-order Taylor estimate on \(B_d(y)\) gives
\begin{equation}\label{eq:display-092}
|Dw_\tau(y)\cdot e|
\le \frac{2}{h}\|w_\tau\|_{C^0(K')}
+\frac{h}{2}\|D^2w_\tau\|_{L^\infty(K')}.
\end{equation}
Taking \(h=\tau^3\) gives 
\begin{equation}\label{eq:display-inline-804}
\|Dw_\tau\|_{C^0(K)}\le C_K\tau^3.
\end{equation}

For \(y,z\in K\) with \(|y-z|<d\), the gradient and Hessian bounds give
\begin{equation}\label{eq:display-093}
|Dw_\tau(y)-Dw_\tau(z)|
\le C_K\min\{\tau^3,|y-z|\}
\le C_K\tau^{3(1-\theta)}|y-z|^\theta .
\end{equation}
For \(|y-z|\ge d\), the same Hölder bound follows from the
supremum bound on \(Dw_\tau\). This proves
\eqref{eq:profilestatement}.

For the gradient asymptotics, take \(K=S^{15}\),
\(\tau=|q-a|\), and \(y=(q-a)/\tau\). Since
\begin{equation}\label{eq:gradient-profile-scaling}
D_yG_\tau=\tau^7D_qG_\Omega,\qquad
\|Dw_\tau\|_{C^0(S^{15})}\le C\tau^3,
\end{equation}
we obtain
\begin{equation}\label{eq:display-094}
D_qG_\Omega(q,a)-6|q-a|^{-8}(q-a)
=\tau^{-7}Dw_\tau(y)=O(\tau^{-4}),
\end{equation}
which proves \eqref{eq:gradientstatement}.
\end{proof}

\subsection*{Acknowledgements}
Xi-Nan Ma was supported by the National Key Research and Development Program of China
(Grant No. 2025YFA1017600) and NSFC No. 12141105.

\subsection*{Statement of usage of artificial intelligence}
The central mathematical ideas and overall approach are due to the authors.
OpenAI ChatGPT was used to assist with proofreading and improving the clarity
of the exposition. The authors take full responsibility for the paper's content
and the accuracy of its references.

\end{document}